\documentclass
{gtpart}
\usepackage{graphicx,hyperref,pinlabel}

\title[Isoperimetric Spectra]{Snowflake groups, Perron--Frobenius
eigenvalues,\\ and isoperimetric spectra} 

\author[N Brady]{Noel Brady}
\givenname{Noel}
\surname{Brady}
\address{Mathematics Department\\ University of Oklahoma\\\newline Norman, OK
73019\\ USA} 
\email{nbrady@math.ou.edu}

\author[M\,R Bridson]{Martin R Bridson}
\givenname{Martin R}
\surname{Bridson}
\address{Mathematical Institute\\ 24-29 St Giles'\\\newline Oxford\\ OX1
3LB\\ UK} 
\email{bridson@maths.ox.ac.uk}

\author[M Forester]{Max Forester}
\givenname{Max}
\surname{Forester}
\address{Mathematics Department\\ University of Oklahoma\\\newline Norman, OK
73019\\ USA}
\email{forester@math.ou.edu}

\author[K Shankar]{Krishnan Shankar}
\givenname{Krishnan}
\surname{Shankar}
\address{Mathematics Department\\ University of Oklahoma\\\newline Norman, OK
73019\\ USA}
\email{shankar@math.ou.edu}


\arxivreference{math.GR/0608155}
\arxivpassword{cmdt8}

\volumenumber{}
\issuenumber{}
\publicationyear{}
\papernumber{}
\startpage{}
\endpage{}
\doi{}
\MR{}
\Zbl{}
\received{}
\revised{}
\accepted{}
\published{}
\publishedonline{}
\proposed{}
\seconded{}
\corresponding{}
\editor{}
\version{}

\makeautorefname{letterthm}{Theorem}
\makeautorefname{lettercor}{Corollary}

\theoremstyle{plain}
\newtheorem{theorem}{Theorem}[section]
\newtheorem{proposition}{Proposition}[section]
\newtheorem{corollary}{Corollary}[section]
\newtheorem{lemma}{Lemma}[section]

\newtheorem{letterthm}{Theorem}

\newtheorem{lettercor}{Corollary}

\makeatletter
  \let\c@lettercor\c@letterthm
  \let\c@theorem\c@subsection
  \let\c@proposition\c@subsection
  \let\c@corollary\c@subsection
  \let\c@lemma\c@subsection
\makeatother

\theoremstyle{definition}
\newtheorem{definition}{Definition}[section]
\newtheorem{remark}{Remark}[section]
\newtheorem{remarks}{Remarks}[section]
\newtheorem*{ack}{Acknowledgments} 
\makeatletter
  \let\c@definition\c@subsection
  \let\c@remark\c@subsection
  \let\c@remarks\c@subsection
\makeatother

\numberwithin{equation}{section}
\makeatletter
  \let\c@equation\c@subsection
  \def\theequation{\bf \thesubsection}
\makeatother

\newcommand{\N}{{\mathbb N}}
\newcommand{\G}{{G}} 
\newcommand{\IP}{{\rm{IP}}} 
\newcommand{\IPk}{{\rm{IP}}^{(k)}} 

\newcommand{\Gs}{{\mathcal G}}
\newcommand{\Ns}{{\mathcal N}}
\newcommand{\abs}[1]{\left| {#1} \right|}
\newcommand{\aabs}[1]{\left\| {#1} \right\|}
\renewcommand{\leq}{\leqslant}
\renewcommand{\geq}{\geqslant}

\renewcommand{\preceq}{\preccurlyeq}
\renewcommand{\succeq}{\succcurlyeq}
\DeclareMathOperator{\area}{Area}
\DeclareMathOperator{\FVol}{FVol}
\DeclareMathOperator{\Vol}{Vol}
\DeclareMathOperator{\id}{id}
\DeclareMathOperator{\shell}{shell}
\def\cprime{$'$}

\begin{document}
 
\begin{abstract} The $k$--dimensional Dehn (or isoperimetric) function of a
group bounds the volume of efficient ball-fillings of $k$--spheres mapped
into $k$--connected spaces on which the group acts properly 
and cocompactly; the bound is given as a function of the volume of the
sphere. We advance significantly the observed range of behavior for such
functions. First, to each non-negative integer matrix $P$ and positive
rational number $r$, we associate a finite, aspherical $2$--complex
$X_{r,P}$ and determine the Dehn function of its fundamental group
$G_{r,P}$ in terms of $r$ and the Perron--Frobenius eigenvalue of 
$P$. The range of functions obtained includes $\delta(x)= x^s$, where
$s\in\Q\cap [2,\infty)$ is arbitrary. Next, special features of the groups
$G_{r,P}$ allow us to construct iterated multiple HNN extensions which
exhibit similar isoperimetric behavior in higher dimensions. In particular, 
for each positive integer $k$ and rational $s\geq (k+1)/k$, there exists a
group with $k$--dimensional Dehn function $x^s$. Similar isoperimetric
inequalities are obtained for fillings modeled on arbitrary manifold
pairs $(M,\partial M)$ in addition to $(B^{k+1},S^k)$. 
\end{abstract}

\maketitle

\section*{Introduction} Given a $k$--connected complex or manifold one
wants to identify functions that bound the volume of efficient
ball-fillings for spheres mapped into that space. The purpose of this
article is to advance  the understanding of which functions
can arise when one seeks optimal bounds in the universal cover of a
compact space. Despite the geometric nature of both the problem and its
solutions, our initial impetus for studying isoperimetric problems  comes
from algebra, more specifically the {\em{word problem}} for groups.

The quest to understand the complexity of word problems has been at the heart
of combinatorial group theory since its inception. When one attacks 
the word problem for a finitely presented group $G$ directly, the most natural
measure of complexity is the Dehn function $\delta(x)$ which bounds 
the number of defining relations that one must apply to a word  $w=_G 1$
to reduce it to the empty word; the bound is a function of word-length
$\abs{w}$. The function $\delta(x)$ recursive if and only if $G$ has a
solvable word problem. 

Progress in the last ten years has led to a fairly complete understanding
of which functions arise as Dehn functions of finitely presented groups. The
most comprehensive information comes from \cite{BRS} where, modulo
issues associated to the $P=NP$ question, Birget, Rips and Sapir essentially
provide  a characterisation of the Dehn functions greater than $x^4$. In 
particular they show that the following {\em isoperimetric spectrum} is
dense in the range $[4,\infty)$.
\[ \IP=\{ \alpha \in [1,\infty) \mid f(x)=x^\alpha \text{ is equivalent
to a Dehn function}\}.\]
Gromov proved that $\IP\cap (1,2)$ is empty and that word hyperbolic
groups can be characterised as those which have linear Dehn functions.
In \cite{BrBr} Brady and Bridson completed the understanding of the {\em
coarse} structure of $\IP$ by providing a dense set of exponents in
$\IP\cap [2,\infty)$. What remains unknown is the fine structure of
$\IP\cap (2,4)$. In particular, it has remained unknown whether $\Q\cap
(2,4)\subset \IP$. There has, however, been recent progress on
understanding Dehn functions below $x^4$ that are not of the form
$x^{\alpha}$.  For instance, Ol{\cprime}shanskii and Sapir \cite{olsap}
have constructed groups with Dehn function $x^2 \log(x)$, and
Ol{\cprime}shanskii \cite{ol} has constructed examples with more exotic,
almost-quadratic behavior.

What Brady and Bridson actually do in \cite{BrBr} is associate to each
pair of positive integers $p\geq q$ a finite aspherical $2$--complex
whose fundamental group $G_{p,q}$ has Dehn function $x^{2\log_2
2p/q}$. These complexes are obtained by attaching a pair of annuli to a
torus, the attaching maps being chosen so as to ensure the existence of
a family of discs in the universal cover that display a certain {\em
snowflake geometry} (cf \fullref{fig:sfdisk} below). In the present
article we present a more sophisticated version of the snowflake
construction that yields a much larger class of isoperimetric exponents. 

\begin{letterthm}\label{mainthm1} 
Let $P$ be an irreducible non-negative integer matrix with
Perron--Frobenius eigenvalue $\lambda>1$, and let $r$ be a rational number
greater than every row sum of $P$. Then there is a finitely presented
group $G_{r,P}$ with Dehn function $\delta(x) \simeq x^{2 \log_{\lambda}
(r)}$. 
\end{letterthm}

Here, $\simeq$ denotes coarse Lipschitz equivalence of functions. 
By taking $P$ to be the $1\times 1$ matrix $(2^{2q})$ and $r=2^{p}$
(for integers $p > 2q$) we obtain the Dehn function $\delta(x) \simeq
x^{p/q}$ and deduce the following corollary. 
\begin{lettercor}\label{maincor1}
$\Q\cap (2,\infty)\subset \IP$. 
\end{lettercor}

The influential work of M.~Gromov \cite{gromov1}, \cite{gromov2} embedded
the word problem in the broader context of filling problems for
Riemannian manifolds and combinatorial complexes.  For example, Gromov's
Filling Theorem \cite{mrb:bfs} states that given a compact Riemannian
manifold $M$, the smallest function bounding the area of least-area
discs in $M$ as a function of their boundary length is coarsely
Lipschitz equivalent to the Dehn function of $\pi_1M$.  In the geometric
context, it is natural to extend questions about the size of optimal
fillings to higher-dimensional spheres, exploring higher-dimensional
isoperimetric functions that bound the volume of optimal ball-fillings
of spheres mapped into the manifold (or complex). Correspondingly, one
defines higher-dimensional Dehn functions $\delta^{(k)}(x)$ for finitely
presented groups $\G$ that have a classifying space with a compact
$(k+1)$--skeleton (see \fullref{prelimsect}). The equivalence class of
$\delta^{(k)}$ is a quasi-isometry invariant of $\G$, by
Alonso--Wang--Pride \cite{AWP}. 

In contrast to the situation of ordinary Dehn functions, Papasoglu
\cite{papa} has shown that $\delta^{(2)}(x)$ is always bounded by a
recursive function. This is not the case in higher dimensions,
however. For each $k > 2$, Young \cite{young} constructs a group for
which $\delta^{(k)}(x)$ is not subrecursive. 

For each positive integer $k$ one has the {\em $k$--dimensional 
isoperimetric spectrum} 
\[ \IPk = \{ \alpha \in [1,\infty) \mid f(x)=x^\alpha \text{ is
equivalent to a $k$--dimensional Dehn function}\}.\] 
We do not yet have as detailed a knowledge of the structure of these sets
as we do of $\IP=\IP^{(1)}$. Indeed knowledge until now has been
remarkably sparse even for $\IP^{(2)}$: the results of Alonso {\em et
al} \cite{alonso+}, Wang and Pride \cite{WP}, and Wang \cite{wang}
provide infinite sets of exponents in the range $[3/2,2)$ and provide
evidence for the existence of exponents in the range $[2,\infty)$; the
snowflake construction of Brady and Bridson \cite{BrBr} provides a dense set
of exponents 
in the interval $[3/2,2)$; and in Bridson \cite{mb:plms} 
it is was proved that $2,3\in\IP^{(2)}$ (see also Burillo
\cite{burillo}). Gromov and others have investigated the isoperimetric
behavior of lattices \cite{gromov2}. 

Our second theorem relieves the dearth  of knowledge about the coarse
structure of $\IP^{(k)}$, $k \geq 2$. 

\begin{letterthm}\label{mainthm2} 
Let $P$ be an irreducible non-negative integer matrix with Perron--Frobenius
eigenvalue $\lambda>1$, and let $r$ be an integer greater than every 
row sum of $P$. Then for every $k \geq 2$ there is a group
$\Sigma^{k-1}G_{r,P}$ of type $\mathcal{F}_{k+1}$ with $k$--dimensional
Dehn function 
$\delta^{(k)} (x) \simeq x^{2 \log_{\lambda}(r)}$. There are also groups
$\Sigma^{k-1} \Z^2$ of type $\mathcal{F}_{k+1}$ with $k$--dimensional Dehn
function $\delta^{(k)}(x) \simeq x^2$. 
\end{letterthm}

By taking $P$ to be the $1\times 1$ matrix $(2^{2q})$ and $r=2^{p}$ we
see that $\Q\cap [2, \infty)\subset \IP^{(k)}$; in particular
 $\IP^{(k)}$ is dense in the range $[2,\infty)$. But
that falls short of one's expectations: as in the case $k=1$, one
anticipates that $\IP^{(k)}$ should be dense in the range that begins
with the exponent $(k+1)/k$ corresponding to the isoperimetric inequality
for spheres in Euclidean space. In order to fulfil this expectation, we
investigate the higher Dehn functions of products $\G\times \Z$ and
prove the following theorem. 

\begin{letterthm}\label{mainthm3} 
If $P,\,\lambda$ and $r$ are as in \fullref{mainthm2}, then for all
$q, \ell\in \N$, the $(q+\ell)$--dimensional Dehn function of
$\Sigma^{q-1}G_{r,P} \times \Z^{\ell}$ is equivalent to $x^{s}$, where $s
= \frac{(\ell + 1)\alpha - \ell}{\ell \alpha - (\ell -1)}$ and $\alpha =
2\log_{\lambda}(r)$. 
The $(q+\ell)$--dimensional Dehn function of $\Sigma^{q-1}\Z^2 \times
\Z^{\ell}$ is equivalent to $x^{s}$, where $s = \frac{\ell+2}{\ell +1}$. 
\end{letterthm}

\begin{figure}[ht!]
\labellist
\pinlabel* {\scriptsize $1$} at 36 5
\pinlabel* {\scriptsize $2$} at 180 5
\pinlabel* {\scriptsize $3$} at 324 5
\pinlabel* {\scriptsize $\frac{3}{2}$} at 108 5
\pinlabel* {\scriptsize $\frac{4}{3}$} at 84 5
\pinlabel* {\scriptsize $\frac{5}{4}$} at 72 5
\pinlabel* {\scriptsize $\frac{6}{5}$} at 64.8 5
\pinlabel* {\scriptsize $\IP$} [l] at 0 22.5
\pinlabel* {\scriptsize $\IP^{(2)}$} [l] at 0 49.5
\pinlabel* {\scriptsize $\IP^{(3)}$} [l] at 0 76.5
\pinlabel* {\scriptsize $\IP^{(4)}$} [l] at 0 103.5
\pinlabel* {\scriptsize $\IP^{(5)}$} [l] at 0 130.5
\pinlabel* {\scriptsize $\nwarrow \times \, \Z$} [bl] at 103.5 57
\pinlabel* {\scriptsize $\uparrow \Sigma$} [bl] at 255 57
\endlabellist
	\begin{center}
		\includegraphics[width=4in]{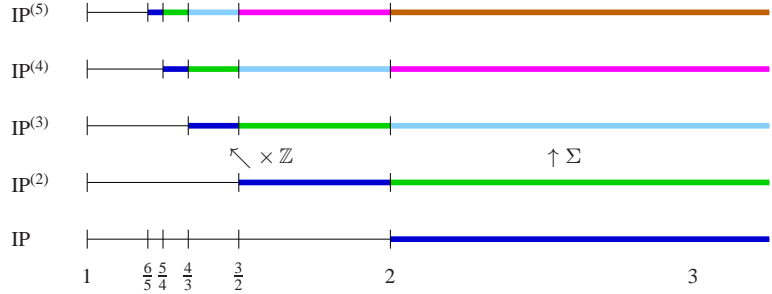}
	\end{center}
	\caption{Isoperimetric exponents of $\Sigma^{q-1} G_{r,P} \times
	\Z^{\ell}$. Colors correspond to fixed values of
$q$.}\label{fig:spectra}  
\end{figure}

By holding $q$ and $\ell$ fixed and varying $r$ and $P$, one obtains
a dense set of exponents $s$ in the interval $[\frac{\ell +2}{\ell+1},
\frac{\ell +1}{\ell}]$ including all rationals in this range. 
By varying $q$ and $\ell$ with $k=q+\ell$ and taking account of
\fullref{mainthm2} we deduce the following result, shown pictorially in 
\fullref{fig:spectra}. 

\begin{lettercor}\label{maincor2} 
$\Q\cap [(k+1)/k, \infty)\subset \IP^{(k)}$. 
\end{lettercor}

The main aim of Brady and Bridson's initial construction of snowflake
groups \cite{BrBr} was to prove that the closure of $\IP^{(1)}$ is
$\{1\}\cup [2,\infty)$. \fullref{maincor2} implies that the closure
of $\IP^{(k)}$ contains $\{1\}\cup [(k+1)/k,\infty)$. Building on this
result, Brady and Forester \cite{bf} have recently shown 
that the closure of $\IP^{(k)}$ is in fact equal to $[1,\infty)$ for $k
\geq 2$. Other examples, known earlier, include solvable and
nilpotent groups whose two-dimensional Dehn functions appear to be $x
\log x$ and $x^{4/3}$ respectively, by Wang \cite{wang} and
Coulhon--Saloff-Coste \cite{csc} (the latter was pointed out to us by
Robert Young). It should be noted, however, that in both cases the upper
bound is derived using the Sobolev inequality from
Varopoulos--Saloff-Coste--Coulhon \cite{vscc}. The resulting
isoperimetric inequality concerns embedded fillings only, which do not
\emph{a priori} suffice for our definitions (which allow singular
maps and fillings). 

This article is organised as follows. In Section 1 we outline the
construction of the snowflake groups $G_{r,P}$ and their HNN extensions 
$\Sigma G_{r,P}$, deferring a detailed account to Sections 4 and 6. In
Section 2 we define the class of maps with which we shall be working and
record some pertinent properties; we also recall those elements of
Perron--Frobenius theory that we require. The groups $G_{r,P}$ are
fundamental groups of graphs of groups; in Section 3 we analyze the
geometry of the vertex groups in these decompositions. The snowflake
geometry of $G_{r,P}$ is described in Section 4 and this is analyzed in
further detail in Section 5 to prove \fullref{mainthm1}. In Section 6
we turn our attention to higher Dehn functions and establish the lower
bounds required for \fullref{mainthm2} by analyzing the geometry of
an explicit sequence of embedded $(k+1)$--balls in the universal cover of
a $(k+1)$--dimensional classifying space for $\Sigma^{k-1}G_{r,P}$. In
Section 7 we establish the complementary upper bounds. The proof proceeds
by induction, slicing balls into slabs based of lower-dimensional
fillings. A lack of control on the topology of these slabs obliges one to
prove a stronger result: instead of establishing bounds only on the
behavior of ball-fillings for spheres, one must establish isoperimetric
inequalities for all pairs of compact manifolds $(M^{(k+1)},\partial M)$
mapping to the space in question. In Section 8 we analyze the
isoperimetric behaviour of products $\G\times\Z$ and complete the
proof of \fullref{mainthm3}. 

\begin{ack}
We would like to thank Josh Barnard for helpful discusssions on
admissible maps, and Robert Young for the observation mentioned above. 
We also thank the referee for many detailed and helpful comments. 

Brady's work was funded in part by NSF grant no.\ DMS-0505707, Bridson's
by a Royal Society Wolfson Research Merit Award and an EPSRC Senior
Fellowship, Forester's by NSF grant no.\ DMS-0605137, and Shankar's by 
NSF grant no.\ DMS-0336681. Brady and Forester thank the Centre de Recerca
Matem\`atica for their support and hospitality. 
\end{ack}

\section{An outline of the basic construction} 

The groups $G_{r,P}$ we consider are fundamental groups of graphs of
groups whose underlying graphs are determined by a non-negative integer
matrix $P$. The edge groups are infinite cyclic, with attaching maps
determined by a rational number $r$. The vertex groups $V_m$ have many
properties in common with free abelian groups of rank $m$. Indeed, for
the purposes of this summary, the reader may take $V_m = \Z^m$ 
(cf~\fullref{vmrmk}). There is a distinguished element $c \in V_m$,
corresponding to the diagonal element (the product of the standard
generators) in $\Z^m$. The precise definition of $V_m$ is given in
\fullref{Vmsect}. 

The key geometric idea behind \fullref{mainthm1} is that efficient
van Kampen diagrams for the groups $G_{r,P}$ exhibit the {\em snowflake
geometry} illustrated in \fullref{fig:sfdisk}. The essential features
of such diagrams are these: the diagram is composed of polygonal
subdiagrams joined across strips so that the dual to the decomposition
is a tree $T$; and each of the polygonal subdiagrams is a van
Kampen diagram in one of the vertex groups $V_m$ (typically it is an 
$(m+1)$--gon with a base labeled by a power of the distinguished $c\in
V_m$ and $m$ other sides labeled by powers of the $m$ standard 
generators of $V_m$). 

The most important class of diagrams are those that are as symmetric as
possible, having the property that as one moves from the circumcenter of
the dual tree to the boundary of the diagram, the joining strips are all
oriented in such a way that the length of the side strip \emph{decreases}
by a factor of $r$ as one journeys towards the boundary. The labels on the
outer sides of the strips are powers of the diagonal elements in various
vertex groups $V_m$, and a crucial feature of our construction is that
the cyclic subgroups $\langle c\rangle\subset G_{r,P}$ are distorted in
a precisely understood manner, with distortion funtion $\simeq x^\alpha$
where $\alpha = \log_\lambda (r)$ and $\lambda$ is the Perron--Frobenius 
eigenvalue of $P$.  This distortion is determined through the 
analysis of certain paths, called {\em snowflake paths}, which play the 
role of quasi-geodesics in $G_{r,P}$. These snowflake paths are the result 
of a curve shortening process; the dynamics of this process are at the
heart of our calculations and this is where the Perron--Frobenius theory
enters -- see \fullref{sfsect}.

If the tree $T$ has radius $d$, then arguing by induction on $d$
in a suitable class of diagrams, one calculates the length of the
boundary to be $\sim d^{k/\alpha}$ if the central polygon has base $\sim
d^k$. One has a precise understanding of the quadratic Dehn functions of
the vertex groups $V_m$, and this leads to an area estimate of $\sim
d^{2k}$ on these diagrams of diameter $\sim d^k$. Thus we obtain a
family of diagrams with area $\sim d^{2k}$ and perimeter $\sim
d^{k/\alpha}$, and an elementary manipulation of logs provides the
required lower bound of $x^{2\log_\lambda(r)}$ for the Dehn
function of $G_{r,P}$. The complementary upper bound is established in
\fullref{thm1sect}.

A key feature in our construction of $G_{r,P}$ is that when $r$ is an
integer, the snowflake diagrams admit a precise scaling by a factor of
$r$, induced by a monomorphism. The ascending HNN extension $G_{r,P}^*$
of $G_{r,P}$ corresponding to this monomorphism is a group in which one
can {\em stack} scaled snowflake diagrams (see \fullref{fig:sfball}). By
putting together two such stacks (using two stable letters) one obtains a
{\em snowflake ball}, having the same proportions as its equatorial
snowflake disk. That is, the interior volume and surface area of the ball
are comparable to the area and boundary length, respectively, of the
equatorial disk. In this way, one discovers the higher dimensional
isoperimetric behavior of the multiple HNN extension $\Sigma G_{r,P}$
which is obtained by amalgamating two copies of $G_{r,P}^*$ along
$G_{r,P}$. 

The snowflake balls just described determine a lower bound for the
$2$--dimensional Dehn function mentioned in \fullref{mainthm2}; upper
bounds in this case can be deduced by using Wang--Pride \cite{WP}. To
proceed in higher dimensions we iterate the suspension procedure
described above. Lower bounds can be determined as before. However, in
dimensions greater than $2$, upper bounds require new techniques. In
particular, we need to consider isoperimetric inequalities for compact
manifolds $(M,\partial M)$ other than $(B^{k},S^{k-1})$. Using this
perspective, we establish general upper bounds for ascending HNN
extensions. This is achieved in \fullref{suspthm} and is further refined
in \fullref{upperGxZ}.

\begin{remark}\label{vmrmk}
The actual vertex groups $V_m$ of $G_{r,P}$ are themselves fundamental 
groups of graphs of groups with vertex groups $\Z^2$ and edge groups
$\Z$. It turns out that this structure is compatible with the larger
$G_{r,P}$ graph of groups structure. That is, $G_{r,P}$ itself may be
viewed as the fundamental group of an aspherical $2$--complex assembled
from a finite collection of tori and annuli. With respect to a fixed
framing on the tori, the attaching maps of the annuli are all powers of
the slopes $\{1/0,\ 0/1,\ 1/1\}$.  From this perspective, it is perhaps
surprising that one can encode the range of isoperimetric exponents
stated in \fullref{mainthm1}.
\end{remark}

\subsection*{An explicit example} \label{example}

We conclude this outine with an explicit example illustrating
\fullref{mainthm1}. The example that we present here has Dehn function
$x^{p/q}$, where $p>2q$ are positive integers (common factors are
allowed). 

Let $P$ be the $1\times 1$ matrix with entry $2^{2q}=4^q$ and let $r=2^p$. 
Then $G_{r,P}$ is the fundamental group of a graph of groups $\mathcal G$
with one vertex group and $4^{q}$ infinite cyclic edge groups.  The
single vertex group $V_{4^{q}}$ is the fundamental group of a
tree of groups that we shall describe in a moment. $V_{4^{q}}$ has
generators $a_1,\dots,a_{4^q}$; the product of these
generators $c=a_1\cdots a_{4^q}$ plays a special role.

The $i$th edge group of $\mathcal G$ has two monomorphisms to the vertex
group $V_{4^{q}}$. One maps the generator to $c$ and the other maps the
generator to $a_i^{2^p}$. Thus we have a relative presentation 
\[ G_{p/q} \ = \ G_{r,P} \ = \ \langle \, V_{4^{q}}, s_{1}, \ldots,
s_{4^{q}} \mid s_{i}^{-1} a_i^{2^{p}}s_{i} = c \ \  (i = 1, \ldots,
4^{q}) \, \rangle. \]

It remains to elucidate the structure of the group $V_{4^{q}}$. This is
the fundamental group of a tree of groups in which each of the vertex
groups is isomorphic to $\Z^2$ and each of the edge groups is
infinite cyclic. The underlying tree is a segment with $4^q - 2$ edges
and $4^q - 1$ vertices. A basis $\{a_i,b_i\}$ is fixed for each vertex
group, and the generator of each edge group maps to the generator $a_i$
of the left-hand vertex group, and to the \emph{diagonal element}
$a_{i+1}b_{i+1}$ of the right-hand vertex group.  

The generators $a_1,\ldots, a_{4^q}$ mentioned above are the
generators $a_i$ of these vertex groups together with $a_{4^q} = b_{4^q -
1}$. The distinguished element $c$ is the diagonal $a_1 b_1$ of the
leftmost vertex group $\Z^{2}$ (see \fullref{fig:relation}(a)). 

\fullref{mainthm1} tells us that the Dehn function of $G_{p/q}$ is
$x^\alpha$ where $ \alpha= 2 \log_{4^q} 2^p = p/q$.  Consider, for example,
the group $G_{5/2}$ with Dehn function $x^{5/2}$. In this case, the tree
described above is a segment of length $14$ and the above description of
$V_{4^q}$ yields the presentation 
\[ \langle a_1, b_1, a_2, b_2, \ldots, a_{15}, b_{15} \mid
[a_i,b_i] \ (i = 1, \ldots, 15), \ b_i = a_{i+1}b_{i+1} \ (i = 1,
\ldots, 14) \, \rangle. \]
Eliminating the superfluous generators $b_1, \ldots, b_{14}$ and
relabelling $b_{15}$ as $a_{16}$, as in the description of $V_{4^q}$
above, we get
\[ V_{16} \ = \ \langle \, a_1,\dots, a_{16} \mid \theta\in\mathcal C_{16}
\, \rangle \]
where $\mathcal C_{16}$ is the following set of commutators: 
\[ [a_1, a_2\cdots a_{16}], \ [a_2, a_3 \cdots a_{16}], \cdots,
[a_{14}, a_{15} a_{16}], \ [a_{15},a_{16}]. \]
Thus we obtain the explicit presentation 
\[ 
G_{5/2} = \langle \, a_1,\ldots, a_{16}, s_1, \ldots, s_{16} \mid \mathcal
C_{16}; \ s_{i}^{-1}a_i^{32}s_{i} = a_1\cdots a_{16} \ (i = 1, \ldots, 16)
\, \rangle. \]

We have just described a 32--generator, 31--relator presentation of
$G_{5/2}$. The corresponding presentation for $G_{p/q}$ has $2^{2q+1}$
generators and $2^{2q+1}-1$ relations.

\section{Preliminaries}\label{prelimsect} 

In the first part of this section we recall the basic definitions
associated to Dehn functions. We then gather those elements
of Perron--Frobenius theory that will be needed in the sequel. 

\subsection*{Dehn functions}
Given a finitely presented group $G=\langle \mathcal A \mid \mathcal R
\rangle$ and a word $w $ in the generators $\mathcal A^{\pm 1}$ that
represents $1\in G$, one defines 
\[ \area (w) \; = \; \min\{ N \in {\N} \mid 
\exists\text{ equality }
w = \prod^N_{j=1}u_jr_j u_j^{-1} \text{ freely,
where } r_j \in \mathcal R^{\pm 1} \}\, .\]
The {\it Dehn function} $\delta(x)$ of the finite presentation $\langle
\mathcal A \mid \mathcal R\rangle$ is defined by 
\[ \delta(x) \; = \; \max\{\text{\rm{Area}}(w) \mid
 w \in \text{ker} (F(\mathcal{A}) \to G), 
\; \abs{w} \leq x \, \} \]
where $\abs{w}$ denotes the length of the word $w$. It is straightforward
to show that the Dehn functions of any two finite presentations of the
same group are equivalent in the following sense (and modulo
this equivalence relation it therefore makes sense to talk of ``the'' Dehn
function of a finitely presented group). 

Given two functions $f,g \co [0, \infty) \to [0, \infty)$ we define $f
\preceq g$ if there exists a positive constant  $C$ such that 
\[ f(x) \leq C\,g(Cx) + Cx \]
for all $x\geq 0$. If $f \preceq g$ and $g \preceq f$ then $f$ and $g$
are said to be \emph{equivalent}, denoted $f \simeq g$. 

\begin{remark} \label{sparse} 
In order to establish the relation $f \preceq g$ between non-decreasing
functions, it suffices to consider relatively sparse sequences of
integers. For if $(n_i)$ is an unbounded sequence of integers for which 
there is a constant $C>0$ such that $n_0 = 0$ and $n_{i+1} \leq Cn_i$ for
all $i$, and if $f(n_i) \leq g(n_i)$ for all $i$, then $f \preceq
g$. Indeed, given $x \in [0, \infty)$ there is an index $i$ such that
$n_i \leq x \leq n_{i+1}$, whence $f(x) \leq f(n_{i+1}) \leq g(n_{i+1})
\leq g(C n_i) \leq g(Cx)$. 
\end{remark}

We refer to Bridson \cite{mrb:bfs} for general facts about Dehn functions, in
particular the interpretation of $\area (w)$ in terms of {\em van Kampen
diagrams} over $\langle \mathcal A \mid \mathcal R \rangle$. Recall that
a van Kampen diagram for $w$ is a labeled, contractible, planar
$2$--complex with a basepoint and boundary label $w$.  Associated to such
a diagram $D$ one has a cellular map $\tilde D$ from $D$ to the universal
cover $\tilde K$ of the standard $2$--complex of $\langle \mathcal A \mid
\mathcal R \rangle$, respecting labels and basepoint. The diagram is said
to be {\em embedded} if this map in injective. 

\begin{remark}\label{embedded} 
If the presentation $\langle \mathcal A \mid \mathcal R \rangle$ is
aspherical and the diagram $D$ is embedded, then $D$ has the smallest
area among all diagrams with the same boundary label. To see this, note
that if $\tilde\Delta$ is a diagram with the same boundary circuit as
$\tilde D$, then $\tilde D - \tilde\Delta$ defines a $2$--cycle in
$\tilde K$, which must be zero since $H_2(\tilde K;\Z) = 0$ and there
are no $3$--cells. Thus  each $2$--cell  in the image of $\tilde D$ must
also occur in the image of $\tilde\Delta$. And since $\tilde D$ is an
embedding, the number of $2$--cells in the image (hence domain) of $\tilde
\Delta$ is at least $\area(D)$. 
\end{remark}

\subsection*{Higher-dimensional Dehn functions}
Our treatment of higher-dimensional Dehn (isoperimetric) functions is
similar to that of Bridson \cite{mb:plms}, which is an interpretation of
the more algebraic treatment of Alonso {\em et al} \cite{AWP}. See
Section 5 of \cite{mb:plms} for an explanation of the differences
with the approaches of other authors, in particular Gromov \cite{gromov2},
Epstein {\em et al} \cite{ECHLPT}, and Hatcher--Vogtmann \cite{HV}. 

The $k$--dimensional Dehn function is a function $\delta^{(k)} \co \N\to \N$
defined for any group $G$ that is of type $\mathcal{F}_{k+1}$ (that 
is, has a $K(G,1)$ with finite $(k+1)$--skeleton). Up to equivalence,
$\delta^{(k)}(x)$ is a quasi-isometry invariant. Roughly 
speaking, $\delta^{(k)}(x)$ measures the number of $(k+1)$--cells that one
needs in order to fill any singular $k$--sphere in $K(G,1)$ comprised of
at most $x$ $k$--cells. The reader who is happy with this description can
skip the technicalities in the remainder of this subsection. However, to
be precise one has to be careful about the classes of maps that one
considers and the way in which one counts cells. To this end, we make 
the following definitions. 

If $W$ is a compact $k$--dimensional manifold and $X$ a CW complex, an
\emph{admissible map} is a continuous map $f \co W \to X^{(k)}\subset X$
such that $f^{-1}(X^{(k)} - X^{(k-1)})$ is a disjoint union of open
$k$--dimensional balls, each mapped by $f$ homeomorphically onto a
$k$--cell of $X$. 

If $f \co W \to X$ is admissible we define the \emph{volume} of $f$,
denoted $\Vol^k(f)$, to be the number of open $k$--balls in $W$ mapping to
$k$--cells of $X$. This notion is useful because of the abundance of
admissible maps: 

\begin{lemma}\label{admissible} 
Let $W$ be a compact manifold (smooth or piecewise-linear) of dimension
$k$ and let $X$ 
be a CW complex. Then every continuous map $f \co W \to X$ is homotopic 
to an admissible map. If $f(\partial W) \subset X^{(k-1)}$ then the homotopy
may be taken rel $\partial W$. 
\end{lemma}

\begin{proof}
We prove the lemma in the smooth case; analogous methods apply in the
piecewise-linear category (cf the transversality theorem of 
Buoncristiano--Rourke--Sanderson \cite{BuRoSa}). 

First arrange that $f(W) \subset X^{(k)}$ using cellular
approximation. Next consider $X^{(k)} - X^{(k-1)}$ as a smooth manifold 
and perturb $f$ to be smooth on the preimage of this open set. Let
$C\subset X^{(k)}$ be a set consisting of one point in the interior of
each $k$--cell of $X$. By Sard's theorem we can choose each point of $C$
to be a regular value of $f$. The preimage $f^{-1}(C)$ is now a 
codimension $k$ submanifold of $W$ (ie a finite set of points) and $f$
is a local diffeomorphism at each of these points, by the inverse
function theorem. Thus there is a neighborhood $V$ of $C$ consisting of
a small open ball around each point, whose preimage in $W$ is a disjoint
union of open balls, each mapping diffeomorphically to a component of
$V$. Now modify $f$ by composing it with a map of $X$ (homotopic to the
identity) that stretches each component of $V$ across the $k$--cell
containing it, and pushes its complement into $X^{(k-1)}$. The resulting
map is admissible. 
\end{proof}

Given a group $G$ of type $\mathcal{F}_{k+1}$, fix an aspherical CW
complex $X$ with 
fundamental group $G$ and finite $(k+1)$--skeleton. Let $\widetilde{X}$ be
the universal cover of $X$. If $f\co S^k \to \widetilde{X}$ is
an admissible map, define the \emph{filling volume} of $f$ to be the
minimal volume of an extension of $f$ to $B^{k+1}$: 
\[ \, \FVol(f) \ = \ \min \{\, \Vol^{k+1}(g) \mid g \co B^{k+1} \to
\widetilde{X}, \  g\vert_{\partial B^{k+1}} = f \, \}.\]  
Note that the maps $g$ must be admissible for volume to be defined. Such
extensions exist by \fullref{admissible}, since $\pi_k(\widetilde{X})$
is trivial. Next we define the \emph{$k$--dimensional Dehn function} of
$X$ to be 
\[\delta^{(k)}(x) \ = \ \sup \{ \, \FVol(f) \mid f \co S^k \to
\widetilde{X}, \ \Vol^k(f) \leq x \, \}.\]
Again, the maps $f$ are assumed to be admissible. We will also write
$\delta^{(k)}(x)$ as $\delta^{(k)}_G(x)$ (recall that $G$ is the
fundamental group of $X$). 

\begin{remarks} \label{dehnremarks} 
(1) In these definitions
one could equally well use $X$ in place of $\widetilde{X}$, since
maps $S^k \to X$ (or $B^{k+1}\to X$) and their lifts to $\widetilde{X}$
have the same volume. There are reasons to prefer $\widetilde{X}$,
however, as we shall see in the next definition below. 

(2) It is not difficult to show that the Dehn function
$\delta_G^{(k)}(x)$ agrees with the notion defined by Alonso {\em et al}
in \cite{AWP}.  A discussion along these lines is given in Bridson
\cite[Section 5]{mb:plms}. Moreover it is proved in \cite{AWP} that, up to 
equivalence, $\delta_G^{(k)}(x)$ depends only on $G$ (and in fact is a
quasi-isometry invariant); hence we refer to it as ``the'' $k$--dimensional
Dehn function of $G$. It is also proved in \cite{AWP} that the supremum
in the definition of $\delta_G^{(k)}(x)$ is attained. 
\end{remarks}

\subsection*{More general Dehn functions} 
The definition of $\delta^{(k)}(x)$ generalizes in a natural way to give
Dehn functions modeled on manifolds other than $B^{k+1}$. For example,
Gromov has defined \emph{genus $g$ filling invariants} based on surfaces
other than the disk \cite{gromov2}. Here we need to consider arbitrary
compact manifolds. 

Let $(M, \partial M)$ be a compact manifold pair (smooth or
piecewise-linear) with
$\dim M = k+1$. If $f \co \partial M \to \widetilde{X}$ is an admissible
map define  
\begin{equation}\label{fvoleqn}
\FVol^M(f) \ = \ \min \{ \, \Vol^{k+1}(g) \mid g \co M \to 
\widetilde{X}, \ g \vert_{\partial M} = f \, \}
\end{equation} 
and
\[ \delta^M(x) \ = \ \sup \{ \, \FVol^M(f) \mid f \co \partial M \to
\widetilde{X}, \ \Vol^k(f) \leq x \, \}.\] 
The \emph{dimension} of $\delta^M(x)$ is $k$, the dimension of $\partial M$
(when $\partial M \not= \emptyset$). In general we do not assume that $M$
is connected or that $\partial M \not= \emptyset$. Note that if $M$ is 
closed then $\delta^M(x)$ is identically zero, since $M$ may be mapped to
a point, of zero volume. We will also use the notation $\delta^M_G(x)$
for $\delta^M(x)$. 

\begin{remarks} 
(1) In the definition of $\delta^M(x)$ it is important that we use maps
into $\widetilde{X}$, which is contractible, since maps $f \co \partial M
\to X$ need not have extensions to $M$. Note that if $(M,\partial M) =
(B^{k+1},S^k)$ then the definitions of $\delta^M(x)$ and $\delta^{(k)}(x)$
agree. 

(2) The omission of $X$ from the notation and the adoption of the
alternative notation $\delta^M_G(x)$ suggest an implicit claim that, as
in the case $M=B^{k+1}$, the equivalence class of $\delta^M(x)$ depends
only on $G$. We shall address this issue elsewhere, as it would take us
too far afield in the context of the current paper. The structure of 
the arguments in Sections 7 and 8 requires us to work with specific
choices of $X$ anyway. 

(3) Also to be addressed elsewhere is whether 
the supremum in the definition of $\delta^M(x)$ is attained. The main
difficulty arises when $M$ is $3$--dimensional, as we shall explain in a
moment. In the current paper this issue plays no role because none of
the bounds that we establish require {\em{a priori}} finiteness.

(4) If $\dim M = k+1 \geq 4$ then $\delta^M(x) \leq \delta^{(k)}(x)$
provided $\partial M$ is connected or $\delta^{(k)}(x)$ is
superadditive.
In particular, $\delta^M(x)$ is finite. 
The key point to observe here is that
 if $N=\partial M$ is connected and  $f:N\to \widetilde X$
has volume $V$,
then there is an admissible homotopy with $(k+1)$--dimensional volume at
most $\delta^{(k)}(V)$
from $f$ to an admissible map
$f':N\to \widetilde X$ whose image lies 
$\widetilde X^{(k-1)}$; one can then fill $f'$
by a map $M\to X$ with zero $(k+1)$--dimensional volume.

To see that this homotopy exists, one considers a
$(k-1)$--sphere $S$ in $N$ that encloses a ball $D$ containing 
all of the open discs that contribute to the volume of $f$.
The restriction of $f$ to $S$ is trivial
in $H_{k-1} (\widetilde X^{(k-1)})$ and hence in
 $\pi_{k-1}(\widetilde X^{(k-1)})$ 
(recall that $\widetilde X^{(k-1)}$ is $(k-2)$--connected, and $k > 2$). 
The null-homotopy $H: B^k\to\widetilde X^{(k-1)}$ of 
$f|_S$ furnished by this observation can be adjoined to $f|_D$ to give
an admissible map $S^k\to \widetilde X$ of volume $V$. This can then
be filled by an admissible map $B^{k+1}\to\widetilde X$ of volume at
most
$\delta^{(k)}(V)$. The desired map $f'$ is defined to be the adjunction
of $f|_{N-D}$ and $H$. 

If $\dim M = 2$ then the same statement holds; this is proved below
in \fullref{2dim}. 
The case $\dim M = 3$ is different: Young \cite{young} has constructed a
group $G$ such that if $M$ is a $3$--manifold with boundary $S^1 \times
S^1$, then $\delta^M(x)$ is strictly larger than
$\delta^{(2)}(x)$. 
\end{remarks}

\begin{remark}\label{k-embedded} 
An obvious adaptation of the argument in \fullref{embedded} shows that
if $X$ is an aspherical $(k+1)$--dimensional CW complex, $g \co M^{k+1}
\to X$ is an embedding, and $f = g \vert_{\partial M}$ (with $f$ and $g$
admissible) then $\FVol^M(f) = \Vol^{k+1}(g)$. That is, the embedding $g$
has minimal volume among all extensions of $f$ to the manifold $M$. We
shall use this fact in particular in the case of high-dimensional balls
to estimate $\delta^{(k)}(x)$ from below. 
\end{remark}

\subsection*{Perron--Frobenius Theory}
A square non-negative matrix $P$ is said to be \emph{irreducible} if for
every $i$ and $j$ there exists $k\geq 1$ such that the $ij$--entry of
$P^k$ is positive. The basic properties of irreducible matrices are 
summarized in the Perron--Frobenius theorem below. See Seneta
\cite{seneta} and Katok--Hasselblatt \cite{KH} for a more thorough
treatment of this theory and its applications. 

\begin{proposition}[Perron--Frobenius theorem] 
Let $P$ be an irreducible non-negative $R \times R$ matrix. Then 
$P$ has one (up to a scalar) eigenvector with positive coordinates
and no other eigenvectors with non-negative coordinates. Moreover, the
corresponding eigenvalue $\lambda$ is simple, positive, and is 
greater than or equal to the absolute value of all other eigenvalues. If
$m$ and $M$ are the smallest and largest row sums of $P$, then $m \leq
\lambda \leq M$, with equality on either side implying equality
throughout. 
\end{proposition}

\begin{lemma}\label{pflemma} 
Let $P$ be an irreducible non-negative $R \times R$ matrix with
Perron--Frobenius eigenvalue $\lambda$. Let $\{v_1, \ldots, v_R\}$ be a
generalized eigenbasis for $P$, with $v_1$ a positive eigenvector for
$\lambda$, and with corresponding inner product $\langle \, \cdot\, , \,
\cdot \, \rangle$. Then $\langle u, v_1 \rangle > 0$ for every
non-negative vector $u \in \R^R - \{0\}$. 
\end{lemma}

\begin{proof}
Decompose $\R^R$ as $W_1 \oplus \cdots \oplus W_k$ where each $W_i$ is a 
generalized eigenspace for $P$, with $W_1 = \langle v_1 \rangle$. Each
$W_i$ is $P$--invariant, as is the non-negative orthant $\Ns$, since 
$P$ is non-negative. The intersection $(W_2 \oplus \cdots \oplus W_k)
\cap \Ns$ must then be trivial, for otherwise it contains an eigenvector
for $P$ other than $v_1$ (or a scalar multiple), by the Brouwer fixed
point theorem. Hence $\langle u, v_1\rangle \not= 0$ for every $u \in \Ns
- \{0\}$. Since $\Ns - \{0\}$ is connected and contains $v_1$, $\langle
u, v_1\rangle$ is positive. 
\end{proof}

\begin{proposition}[Growth rate]\label{pfprop} 
Let $P$ be an irreducible non-negative $R \times R$ matrix with
Perron--Frobenius eigenvalue $\lambda$. Let $\aabs{ \, \cdot \, }$ be a
norm on $\R^R$. Then there are positive
constants $A_0, A_1$ such that for every non-negative vector $u$ in
$\R^R$ and every integer $k>0$, \ $A_0 \lambda^k\aabs{u} \ \leq \ \aabs{P^k
u} \ \leq \ A_1 \lambda^k\aabs{u}$. 
\end{proposition}

\begin{proof}
First, it is clear that by varying the constants, it suffices to consider
any single norm $\aabs{ \, \cdot \, }$. Consider a generalized eigenbasis
$\{v_1, \ldots, v_R\}$ as in \fullref{pflemma} (with $v_1$ a positive
eigenvector for $\lambda$). Let $\langle \, \cdot \, , \, \cdot \,
\rangle$ and $\aabs{ \, \cdot \, }$ be the corresponding inner product
and norm on $\R^R$. Let $\pi \co \R^R \to \langle v_1 \rangle$ be
orthogonal projection ($\pi(u) = \langle u, v_1 \rangle v_1$). 

Define $A_0 = \inf\{ \aabs{\pi(u)}/\aabs{u} \mid u \in \Ns - \{0\}\}$. Note
that $A_0 > 0$ by \fullref{pflemma} and compactness of $\Ns - \{0\}$
modulo homothety. For every $u \in \Ns - \{0\}$ we now have $\lambda^k A_0
\aabs{u} \leq \lambda^k \aabs{\pi(u)} = \aabs{P^k \pi (u)} \leq \aabs{P^k
u}$. We also have $\aabs{P^k u} \leq \lambda^k \aabs{u}$ since $\lambda$
is the spectral radius of $P$; hence $A_1 = 1$ will work. 
\end{proof}

\section{The vertex groups $V_m$}\label{Vmsect} 

In this section we define groups $V_m$ for each integer $m \geq 2$.  We
begin with a very brief overview of the construction of the groups
$G_{r,P}$ so that the reader knows where the groups $V_{m}$ fit into the
overall picture. 

An irreducible matrix $P$ determines a directed graph (whose 
transition matrix is $P$). This graph is the underlying graph in a 
graph of groups description of the $G_{r, P}$ in \fullref{mainthm1}. 
The vertex groups in this graph of groups are precisely the groups
$V_{m}$ which we define and study in this section. 

The groups $V_{m}$ satisfy a number of the properties that the free
abelian groups $\Z^{m}$ do, but they have geometric dimension 2. In
particular, $V_m$ has generators $a_1, \ldots, a_m$ and has the
following scaling property (cf equation \eqref{equations}): for any
integer $N>0$, the equality $a_1^N \cdots a_m^N = (a_1 \cdots a_m)^N$ 
holds. Moreover, this equality requires on the order of $N^{2}$ relations
of $V_{m}$.  This follows as a special case of \fullref{Vmarea}, which
gives careful estimates on the areas of certain words in $V_{m}$. 

\subsection*{The groups $V_m$} 
\begin{figure}[ht!]
\begin{center}
\labellist
\hair 2pt
\pinlabel {\small $c$} [b] at 95 21
\pinlabel {\small $a_1$} [br] at 13 53 
\pinlabel {\small $a_2$} [br] at 65 101 
\pinlabel {\small $a_3$} [bl] at 131 101 
\pinlabel {\small $a_4$} [Bl] at 180 57 
\pinlabel {\small $b_1$} [tr] at 105 52 
\pinlabel {\small $b_2$} [r] at 136 65
\pinlabel {\small (a)} [t] at 97 9 
\pinlabel {\small (b)} [t] at 304 9
\endlabellist
\includegraphics[width=5in]{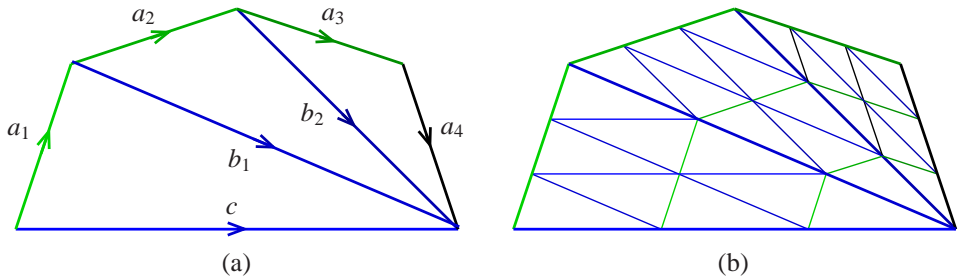}
\end{center}
\caption{Some relations in $V_4$: $c = a_1 a_2 a_3 a_4$ and $c^3 = (a_1)^3
(a_2)^3 (a_3)^3 (a_4)^3$} \label{fig:relation} 
\end{figure}

Begin with $m-1$ copies of $\Z \times \Z$, the $i$th copy having
generators $\{a_i, b_i\}$. The group $V_m$ is formed by successively
amalgamating these groups along infinite cyclic subgroups by adding 
the relations 
$$ b_1 = a_2 b_2, \quad b_2 = a_3 b_3, \quad \ldots, \quad b_{m-2} = a_{m-1}
b_{m-1}.$$ 
Thus $V_m$ is the fundamental group of a graph of groups whose underlying
graph is a segment having $m-2$ edges and $m-1$ vertices. We define two new
elements: $c = a_1 b_1$ and $a_m = b_{m-1}$. Then $a_1, \ldots, a_m$
generate $V_m$ and the relation $a_1 \cdots a_m=c$ holds; see
\fullref{fig:relation}(a). The element $c$ is called the \emph{diagonal 
element} of $V_m$. The additional relations $b_{m-2} = a_{m-1}a_m,
\ \ldots, \ b_{m-k} = a_{m-k+1} \cdots a_m$ are also evident from
\fullref{fig:relation}(a). 

If $m=1$ then we define $V_m$ to be the infinite cyclic group $\langle
a_1 \rangle$ and we set $c = a_1$. Lemmas \ref{shuffle} and
\ref{Vmarea} below clearly hold in this case. 

\begin{lemma}[Shuffling Lemma]\label{shuffle} 
Let $w = w(a_1, \ldots, a_m, c)$ be a word representing $c^N$ in $V_m$
for some integer $N$. Let $n_i$ be the exponent sum of $a_i$ in $w$, and
$n_c$ the exponent sum of $c$ in $w$. Then the words $a_1^{n_1} \cdots 
a_m^{n_m}c^{n_c}$ and $c^{n_c}a_m^{n_m} \cdots a_1^{n_1}$ also represent
$c^N$ in $V_m$ and $n_i = N - n_c$ for all $i$. 
\end{lemma}

\begin{proof}
First we prove the second statement. The abelianization $V_m/[V_m,V_m]
\cong \Z^m$ has $\{a_1, \ldots, a_m\}$ as a basis and the image of $w$
is $a_1^{n_1 + n_c} \cdots a_m^{n_m + n_c}$. Since $c^N$ abelianizes to
$a_1^N \cdots a_m^N$, we must have $n_i = N - n_c$ for all $i$. 

To prove the first statement it now suffices to establish the following
set of equalities for any integer $N$: 
\begin{equation}\label{equations} 
(a_1 \cdots a_m)^N \ = \ a_1^N\cdots a_m^N \ = \ a_m^N \cdots a_1^N \ =
  \ (a_m\cdots a_1)^N.
\end{equation} 
In fact we shall prove the following equalities, by induction on $k$: 
\begin{equation*}\label{shuffleeqn} 
(a_{m-k+1} \cdots a_m)^N \ = \ {a_{m-k+1}}^N\cdots a_m^N \ = \ a_m^N \cdots
  {a_{m-k+1}}^N \ = \ (a_m\cdots a_{m-k+1})^N.
\end{equation*}
The case $k=1$ is evidently true. Suppose the equations hold for a given
$k\geq 1$. By the induction hypothesis $a_{m-k}^N a_{m-k+1}^N \cdots
a_m^N = a_{m-k}^N (a_{m-k+1} \cdots a_m)^N$. Then since $b_{m-k} =
a_{m-k+1} \cdots a_m$ and this element commutes with $a_{m-k}$, we 
conclude that 
$a_{m-k}^N (a_{m-k+1} \cdots a_m)^N = (a_{m-k} \cdots a_m)^N$. The same
commutation relation also yields 
\begin{equation*}
\begin{split}
a_{m-k}^N (a_{m-k+1} \cdots a_m)^N \ &= \ (a_{m-k+1} \cdots a_m)^N
a_{m-k}^N \\
&= \ (a_m \cdots a_{m-k+1})^N a_{m-k}^N \\
&= \ a_m^N \cdots a_{m-k+1}^N a_m^N. 
\end{split}
\end{equation*}
Finally we have $(a_m \cdots a_{m-k+1})^N a_{m-k}^N = (a_m \cdots a_{m-k+1}
a_{m-k})^N$, again because $a_{m-k}$ and $b_{m-k}$ ($=a_m \cdots
a_{m-k+1}$) commute. 
\end{proof}

\begin{remark}[Scaling in $V_m$]\label{scaling}
Equation \eqref{equations} plays a key role in this article. It shows that
the basic relation shown in \fullref{fig:relation}(a) holds 
at larger scales as well. \fullref{fig:relation}(b) illustrates how
these larger relations follow from the triangular relations $b_{i-1} =
a_i b_i$ and $b_{i-1} = b_i a_i$. 
\end{remark}

\subsection*{The spaces $X_m$} 
To compute area in $V_m$ we shall use a specific aspherical $2$--complex
$X_m$ with fundamental group $V_m$. This complex is a union of $m-1$
tori, each triangulated with two $2$--cells realizing the relations
$a_ib_i = b_{i-1}$ and $b_i a_i = b_{i-1}$ (where $b_0 = c$ in the case
$i=1$). Thus the $i$th torus has standard generators given by the
$1$--cells $a_i$ and $b_i$, and its diagonal is joined to the $1$--cell
$b_{i-1}$ of the previous torus. In all there is one vertex, $1$--cells
$a_1, \ldots, a_{m-1}, b_0, \ldots, b_{m-1}$, and $2(m-1)$ triangular
$2$--cells. 

The universal cover $\widetilde{X}_m$ is a union of planes, each covering
one of the tori below. Each plane contains three families of parallel
lines covering the $1$--cells $a_i$, $b_i$, and $b_{i-1}$. The plane
intersects neighboring planes along the $b_j$--lines for $j \not= 0,
m-1$. These planes are the \emph{vertex spaces} of $\widetilde{X}_m$
corresponding to the graph of groups decomposition of $V_m$ described
earlier. The incidence graph of the vertex spaces is the Bass--Serre tree
for this decomposition, with edges corresponding to $b_j$--lines ($j \not=
0, m-1$). 

\begin{remark}\label{Vmdisks} 
\fullref{fig:relation}(b) shows an embedded disk in $\widetilde{X}_m$
with boundary word of the form $c^N = a_1^N \cdots a_m^N \ (N=3)$. The
triangles shown are $2$--cells of $\widetilde{X}_m$. Each large triangular
region lies in a vertex space of $\widetilde{X}_m$. There are similar
embedded disks with boundary word $c^N = a_m^N \cdots a_1^N$ as well. All
of these disks have area $(m-1)N^2$. 
\end{remark}

Throughout this article we usually work with the \emph{standard generators}
$\{a_1, \ldots, a_m\}$ for $V_m$. However in the area computation
below we allow words involving the elements $b_i$ as well. 

\begin{lemma}[Area in $V_m$]\label{Vmarea}
Let $w(a_1, \ldots, a_{m-1}, b_0, \ldots, b_{m-1})$ be a word
representing the element $x^N$ for some $N$, where $x$ is a generator 
$a_i$ or $b_i$. Let $w$ be expressed as $w_1 \cdots w_k$ where each $w_i$
is a power of a generator. Then $N \leq \abs{w}$ and $\area(w x^{-N})
\leq 3 \sum_{i<j} \abs{w_i}\abs{w_j}$. 
\end{lemma}

Note that if the sum included diagonal terms of the form $(3/2) 
\abs{w_i}^2$ then the area bound would simply be $(3/2)\abs{w}^2$. The
leeway afforded by the absence of these terms will be exploited
in the proof of \fullref{mainthm1}. (In particular, it would not
suffice to know only that $V_m$ has quadratic Dehn function.) Also the
statement $N \leq \abs{w}$ implies that every vertex space is a totally
geodesic subspace of $\widetilde{X}_m$. 

\begin{proof}
First we prove that $N \leq \abs{w}$ and then we establish the area
bound. Both proofs are by induction on the complexity of the word $w$, 
defined as follows. Let $p$ be a path in the $1$--skeleton of
$\widetilde{X}_m$ whose edge labels read $w$. Since $w$ represents $x^N$,
the endpoints of $p$ lie in a single vertex space. Hence the induced path
$\widehat{p}$ in the Bass--Serre tree is a closed path. The
\emph{complexity} of $w$ is the length of $\widehat{p}$. Note
that vertices of $\widehat{p}$ correspond to edges of $p$ (or letters of
$w$) and edges correspond to transitions between certain pairs of
generators. Thus the complexity is also the number of such transitions
occurring in $w$. 

If $w$ has complexity zero then $p$ lies in a plane. The statement $N
\leq \abs{w}$ amounts to saying that $x^N$ is a geodesic, which is
clear. If $\widehat{p}$ has positive length then there is a non-trivial
proper subpath $p' \subset p$ with endpoints on a single $b_j$--line. 
(These endpoints correspond to edges in $\widehat{p}$ that map to the same
edge of the Bass--Serre tree, crossing and returning.) 
The subword $w'\subset w$ corresponding to
$p'$ represents an element of the form ${b_j}^M$. Let $u$ be the word
obtained from $w$ by substituting ${b_j}^M$ for $w'$. Then $u$ and $w'$
both have complexity strictly smaller than that of $w$.  By the induction
hypothesis, $M \leq \abs{w'}$ and $N \leq \abs{u} = (\abs{w} - \abs{w'})
+ M$. Therefore $N \leq \abs{w}$. 

Next we establish the area bound when $w$ has complexity zero. Since $p$
then lies entirely within a vertex space of $\widetilde{X}_m$, we may
assume without loss of generality that $V_m = V_1$ and $x = b_0$, so that
$w(a_1, b_0, b_1) = {b_0}^N$ in $V_1 = \langle a_1, b_1, b_0 \mid a_1 b_1
= b_0 = b_1 a_1 \rangle$. Since this group is abelian we can successively
transpose adjacent subwords $w_i$ and cancel pairs of the form $x
x^{-1}$, to obtain $v = a_1^n b_1^n b_0^{N-n}$ for some $n$. Each
transposition of letters contributes $2$ to $\area(w v^{-1})$, so we have
$\area (w v^{-1}) \ \leq \ 2\sum_{i<j} \abs{w_i} \abs{w_j}$. Next let
$I_a$ and $I_b$ be the sets of indices for which $w_i$ is a power
of $a_1$ and $b_1$ respectively. Then $\sum_{i \in I_a} \abs{w_i} \geq
\abs{n}$ and $\sum_{i \in I_b} \abs{w_i} \geq \abs{n}$, and therefore
$\sum_{i<j} \abs{w_i} \abs{w_j} \geq n^2 = \area(vb_0^{-N})$. Then we
have $\area(w b_0^{-N}) \leq \area(w v^{-1}) + \area(v b_0^{-N}) \leq
3\sum_{i<j} \abs{w_i} \abs{w_j}$ as desired. 

Now suppose $w$ has positive complexity. Define $w' \subset w$ and $u$ as
before, so that $w'$ represents ${b_j}^M$, $u$ is obtained from $w$ by
substituting ${b_j}^M$ for $w'$, and both $u$ and $w'$ have smaller
complexity than $w$. Note that $w' = w_{i_0} \cdots w_{i_1} \subset w_1
\cdots w_k$ for some $i_0$ and $i_1$, and so $u = w_1 \cdots w_{i_0 - 1}
{b_j}^M w_{i_1 +1} \cdots w_k$. Let $I = \{i_0, \ldots, i_1\}$. Applying
the induction hypothesis to $u$ and $w'$ we obtain 
\begin{equation}\label{uarea} 
\area(ux^{-N}) \ \leq \ 3 \sum_{\substack{i<j \\ i,j \not\in I}}
\abs{w_i} \abs{w_j} \ + \ 3 \sum_{i \not\in I} \abs{w_i} M 
\end{equation}
and
\begin{equation}\label{w'area} 
\area(w' {b_j}^{-M}) \ \leq \ 3 \sum_{\substack{i<j \\ i,j \in I}} \abs{w_i}
\abs{w_j}. 
\end{equation}
Since $M \leq \abs{w'} = \sum_{j \in I}\abs{w_j}$, inequality 
\eqref{uarea} becomes
\begin{equation}\label{uarea2} 
\area(ux^{-N}) \ \leq \ 3 \sum_{\substack{i<j \\ i,j \not\in I}}
\abs{w_i} \abs{w_j} \ + \ 3 
\Bigl(\sum_{i\not\in I} \abs{w_i} \Bigr)
\Bigl(\sum_{j \in I} \abs{w_j}\Bigr). 
\end{equation}
Adding together \eqref{w'area} and \eqref{uarea2} yields 
\begin{equation*}
\area(w'{b_j}^{-M}) + \area(u x^{-N}) \ \leq \ 3 \sum_{i<j} \abs{w_i}
\abs{w_j} 
\end{equation*}
which proves the lemma because $\area(wx^{-N}) \leq \area(w u^{-1}) +
\area(u x^{-N})$ and $\area(w u^{-1}) = \area(w' {b_j}^{-M})$. 
\end{proof}

\section{The groups $G_{r,P}$ and snowflake words}\label{sfsect}

\subsection*{The groups $G_{r,P}$} 
Start with a non-negative square integer matrix $P = (p_{ij})$ with $R$
rows. Let $m_i$ be the sum of the entries in the $i$th row and let $n =
\sum_i m_i$, the sum of all entries. Form a directed graph $\Gamma$ with
vertices $\{v_1, \ldots, v_R\}$ and having $p_{ij}$ directed edges from
$v_i$ to $v_j$. Label the edges as $\{e_1, \ldots, e_n\}$ and define two
functions $\rho, \sigma \co \{1, \ldots, n\} \to \{1, \ldots, R\}$
indicating the initial and terminal vertices of the edges, so that $e_i$
is a directed edge from $v_{\rho(i)}$ to $v_{\sigma(i)}$ for each
$i$. These functions also indicate the row and column of the matrix entry
accounting for $e_i$. Partition the set $\{1, \ldots, n\}$ as $\bigcup_i
I_i$ by setting $I_i = \rho^{-1}(i)$. Note that $\abs{I_i} = m_i$. 

Let $M = \max \{m_i\}$ and choose a rational number $r = p/q$ with $p > M q
> 0$. We define a graph of groups $\Gs_{r,P}$ with underlying graph
$\Gamma$ as follows. The vertex group $G_{v_i}$ at $v_i$ will be
$V_{m_i}$ and all edge groups will be infinite cyclic. Relabel the
standard generators of these vertex groups as $\{a_1, \ldots, a_n\}$ in
such a way that the standard generating set for $G_{v_i}$ is $\{ a_j \mid
j \in I_i\}$. Let $c_i$ be the diagonal element of the vertex group
$G_{v_i}$. Then the inclusion maps are defined by mapping the generator
of the infinite cyclic group $G_{e_i}$ to the elements ${a_i}^p \in
G_{v_{\rho(i)}}$ and ${c_{\sigma(i)}}^q \in G_{v_{\sigma(i)}}$. 

Let $s_i$ be the stable letter associated to the edge $e_i$. The
fundamental group $G_{r,P}$ of $\Gs_{r,P}$ is obtained from the
presentation 
\begin{equation*}
\langle \, G_{v_1}, \ldots, G_{v_R}, s_1, \ldots, s_n \mid s_i^{-1}
{a_i}^p s_i = {c_{\sigma(i)}}^q \text{ for all } i\, \rangle 
\end{equation*} 
 by adding relations $s_i = 1$ for each edge $e_i$ in a maximal tree in
$\Gamma$. However, we shall continue to use the generating set $\{a_1, 
\ldots, a_n, s_1, \ldots, s_n\}$ for $G_{r,P}$ even though some of these
generators are trivial. 

\subsection*{The spaces $X_{r,P}$} 
We define aspherical $2$--complexes $X_{r,P}$ by forming graphs of spaces
modeling $\Gs_{r,P}$. Namely, take the disjoint union of the spaces
$X_{v_i} \approx X_{m_i}$ (one for each vertex $v_i$) and attach annuli
$A_i$, one for each edge $e_i$ of the graph. The two boundary curves of
$A_i$ are attached to the paths labeled ${a_i}^p$ in $X_{v_{\rho(i)}}$
and ${c_{\sigma(i)}}^q$ in $X_{v_{\sigma(i)}}$. The resulting $2$--complex
$X_{r,P}$ has fundamental group $G_{r,P}$ and it is aspherical because it
is the total space of a graph of aspherical spaces. 

The universal cover $\widetilde{X}_{r,P}$ is a union of copies of the
universal covers $\widetilde{X}_{v_i}$ and infinite strips $\R \times
[-1, 1]$ covering the annuli $A_i$. Each strip is tiled by $2$--cells
whose boundary labels read $s_i^{-1} {a_i}^p s_i {c_{\sigma(i)}}^{-q}$;
the two sides $\R \times \{\pm 1\}$ consist of edges labeled $a_i$ and
$c_{\sigma(i)}$ respectively. Note that if a path crosses a strip along
an edge labeled $s_i$ and returns over $s_i^{-1}$ then the power of $a_i$
represented by the path is divisible by $p$. 

\subsection*{Snowflake words} 
For each group element of the form $c_i^N$ we will define two types of
words in the generators $\{a_1, \ldots, a_n, s_1, \ldots, s_n\}$
representing that element, called positive and negative snowflake
words. The structure of these words is governed by the dynamics of the
matrix $P$. Some snowflake words are close to geodesics, and these are
useful in determining the large scale geometry of $G_{r,P}$. 

We define snowflake words recursively on $\abs{N} \in \N$ as follows. Let

\[N_0 \ = \ \frac{p(M(q+2 + p))}{p-Mq} + p.\] 
Note for future reference that $N_0 > p$. Let
$c$ be the diagonal element of a vertex group with standard ordered
generating set $\{a_{i_1}, \ldots, a_{i_m}\}$. A word $w$ representing
$c^N$ is a \emph{positive snowflake word} if either 
\begin{enumerate}
\item[\textup{(i)}] $\abs{N} \leq N_0$ and $w = a_{i_1}^N \cdots
a_{i_m}^N$, or 
\item[\textup{(ii)}] $\abs{N} > N_0$ and $w = (s_{i_1} u_1 s_{i_1}^{-1})
(a_{i_1}^{N_1}) \cdots (s_{i_m} u_m s_{i_m}^{-1}) (a_{i_m}^{N_m})$ where
each $u_j$ is a positive snowflake word representing a power of
$c_{\sigma(i_j)}$ and $\abs{N_j}<p$ for all $j$. 
\end{enumerate}
In the second case note that each subword $(s_{i_j} u_j
s_{i_j}^{-1})(a_{i_j}^{N_j})$ represents a power of $a_{i_j}$, and by
\fullref{shuffle} this power is $N$. Then since $\abs{N_j} < p$, the
word $(s_{i_j} u_j s_{i_j}^{-1})$ represents either ${a_{i_j}}^{\lfloor
N/p\rfloor p}$ or ${a_{i_j}}^{\lceil N/p \rceil p}$. Consequently, the
word $u_j$ represents either ${c_{\sigma(i_j)}}^{\lfloor N/p\rfloor q}$
or ${c_{\sigma(i_j)}}^{\lceil N/p \rceil q}$. Recall that $\lfloor x 
\rfloor$ and $\lceil x \rceil$ denote the integers closest to $x$ such
that $\lfloor x \rfloor \leq x \leq \lceil x \rceil$, and so $\lfloor
N/p\rfloor p$ and $\lceil N/p \rceil p$ are the multiples of $p$ nearest
to $N$. 

A \emph{negative snowflake word} is defined similarly, with the
ordering of the terms representing powers of $a_{i_j}$ reversed. More
specifically, $w$ satisfies either 
\begin{enumerate}
\item[\textup{(i$'$)}] $\abs{N} \leq N_0$ and $w = a_{i_m}^N \cdots
a_{i_1}^N$, or 
\item[\textup{(ii$'$)}] $\abs{N} > N_0$ and $w = (a_{i_m}^{N_m})(s_{i_m}
u_m s_{i_m}^{-1})  \cdots (a_{i_1}^{N_1})(s_{i_1} u_1 s_{i_1}^{-1})$
where $u_j$ is a negative snowflake word representing a power of
$c_{\sigma(i_j)}$ and $\abs{N_j}<p$ for all $j$. 
\end{enumerate}
As with positive snowflake words, each word $u_j$ will represent either
${c_{\sigma(i_j)}}^{\lfloor N/p\rfloor q}$ or ${c_{\sigma(i_j)}}^{\lceil
N/p \rceil q}$. 

\begin{figure}[ht!]
\begin{center}
\labellist
\hair 2pt
\pinlabel {\small ${c_{\sigma(1)}}^{\lfloor N/p\rfloor q}$} [r] at 70 83
\pinlabel {\small ${c_{\sigma(2)}}^{\lfloor N/p\rfloor q}$} [b] at 189 159
\pinlabel {\small ${c_{\sigma(3)}}^{\lceil N/p\rceil q}$} [bl] at 321 73
\pinlabel {\small $a_1^N$} [tl] at 95 91
\pinlabel {\small $a_2^N$} [t] at 191 136
\pinlabel {\small $a_3^N$} [r] at 290 81
\pinlabel {\small $a_1^{N_1}$} [r] at 81 146
\pinlabel {\small $a_2^{N_2}$} [l] at 277 169
\pinlabel {\small $a_3^{N_3}$} [l] at 363 22
\pinlabel {\small $c^N$} [b] at 196 21
\endlabellist
\includegraphics[width=3.5in]{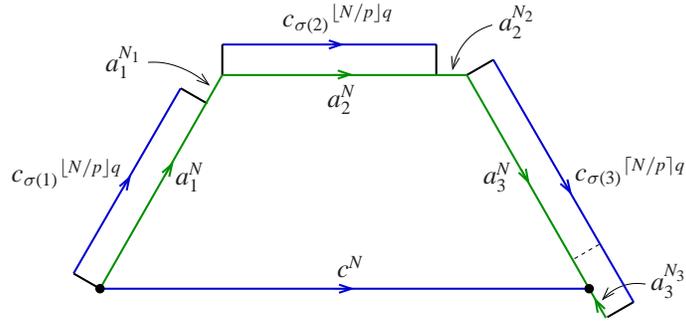}
\end{center}
\caption{One way of shortening $c^N$. Here $\{a_1, a_2, a_3\}$ is the
generating set for a vertex group $V_3$ with diagonal element $c$. The
exponents $N_1$ and $N_2$ are both $N - \lfloor N/p \rfloor p$ and $N_3$
is $N - \lceil N/p \rceil p$. The short black edges are labeled $s_1$,
$s_2$, $s_3$.} \label{fig:curveshortening} 
\end{figure}

To see that the recursion is well-founded note that the definition
describes an iterated curve shortening process in which subwords of the
form $c^N$ are replaced by the words described in case \textup{(ii)} or
\textup{(ii$'$)}, with appropriate powers of $c_{\sigma(i_j)}$ in place 
of $u_j$; see \fullref{fig:curveshortening}. 
Writing $\abs{N} = Ap + B$ with $0 \leq B < p$, the new word representing
$c^N$ has length at most 
\[M \cdot\max\{Aq + 2 + B, (A+1)q + 2 + (p-B)\} \ \leq \ M((A+1)q + 2 + p).\] 
The latter quantity is stictly less than $\abs{N} = Ap + B$ provided $A
(p-Mq) > M(q+2+p)$. Since $B < p$, this occurs whenever $\abs{N} \geq
N_0$. Thus, the new curve is strictly shorter than $c^N$ if $\abs{N} \geq
N_0$. Eventually the subwords $c^N$ all
have length at most $N_0$ and the shortening procedure terminates. See
also \fullref{fig:sfdisk} for the end result of this process. In this
figure the top and bottom halves of the boundary are positive and
negative snowflake words representing $c^N$. 

Note that every snowflake word has a nested structure in which various
subwords are themselves snowflake words. These are the subwords $u_j$
arising at each stage. The minimal such subwords are those given by
\textup{(i)} and \textup{(i$'$)} and these will be called \emph{terminal
subwords}. The \emph{depth} of a snowflake subword is the number of
snowflake subwords of type \textup{(ii)} or \textup{(ii$'$)} properly
containing it, including the original snowflake word
itself. Equivalently, it is the number of matching $s_j$, $s_j^{-1}$
pairs enclosing it. Note that a snowflake word $w$ contains a depth zero
terminal subword if and only if $w$ has the form \textup{(i)} or
\textup{(i$'$)}. 

It is worth emphasizing that the curve shortening process is not
canonically determined, but allows many choices. In each ``remainder''
term $a_i^{N_i}$ the exponent $N_i$ may be positive or negative; the two
possible values for $N_i$ are $N - \lfloor N/p\rfloor p$ and $N - \lceil
N/p \rceil p$. \fullref{fig:curveshortening} shows both possibilities
occurring in a single step, for example.  For this reason, a single
snowflake word may have terminal subwords of different depths. However,
\fullref{sfdepth} below shows that these depths will not differ
substantially. 

\begin{remark}\label{integerpower} 
A special type of snowflake word plays a key role in 
the proof of \fullref{mainthm2}. If $r$ is an 
integer (that is, $r = p/1$) and $N = r^k$ for some $k$, then the
positive (resp. negative) snowflake word representing $c_i^N$ is
unique. What happens is that the exponents $N_j$ in the expressions
\textup{(ii)} or \textup{(ii$'$)} at each stage are always zero; there
are no ``remainder'' terms $a_{i_j}^{N_j}$. Each subword $u_j$ represents
${c_{\sigma(i_j)}}^{N/r}$, and $N/r$ is again a power of $r$. Furthermore,
all terminal subwords will have the form $a_{i_1} \cdots a_{i_m}$ or
$a_{i_m} \cdots a_{i_1}$. 
\end{remark}

\begin{lemma}[Terminal subword depth] \label{sfdepth} 
Given $r$ and $P$ there are positive constants $B_0, B_1$ with the
following property. If a non-trivial snowflake word $w$ representing 
$c^N$ contains a terminal subword of depth $d$ then $B_0 r^d \leq \abs{N}
\leq B_1 r^d$. 
\end{lemma}

\begin{proof}
If $d = 0$ then $w$ has the form \textup{(i)} or \textup{(i$'$)} and $1
\leq \abs{N} \leq N_0$. Thus we need to arrange that $B_0 \leq 1$ and $B_1
\geq N_0$ for the lemma to hold in this case. 

If $d > 0$ then we will show by induction on $d$ that 
\begin{equation} \label{depthbound} 
N_0 r^{d-1} - p (r^{d-2} + \cdots + r + 1) \ \leq \ \abs{N}
\ \leq \ N_0 r^d + p(r^{d-1} + \cdots + r + 1).
\end{equation}
The lower bound then gives 
\[ \abs{N} \ \geq \ N_0 r^{d-1} - p\left(\frac{r^{d-1} - 1}{r-1}\right)
\ \geq \ \frac{1}{r}\left( N_0 - \frac{p}{r-1}\right) r^{d}.\]
Recall that $N_0 > p$ and $r \geq 2$, which imply $N_0 > p/(r-1)$. Now
we may find $B_0 > 0$ so that $B_0 \leq r^{-1}(N_0 - p/(r-1))$ and $B_0
\leq 1$, giving the desired bound. 

The upper bound in \eqref{depthbound} gives
\[ \abs{N} \ \leq \ N_0 r^d + p \left(\frac{r^d - 1}{r-1}\right) \ \leq
\ (N_0 + p) r^d\] 
where the last inequality uses the fact that $r-1 \geq 1$. Now choose
$B_1 \geq N_0 + p$ to obtain the desired bound. 

Next we prove \eqref{depthbound} by induction on $d$. If $d=1$ then
$\abs{N} > N_0$ and $w$ is of the form \textup{(ii)} or \textup{(ii$'$)} 
where some $u_j$ has the form \textup{(i)} or \textup{(i$'$)}. Then $u_j$
represents ${c_{\sigma(i_j)}}^{N'}$ with $N' \leq N_0$, and so $(s_{i_j}
u_j s_{i_j}^{-1})$ represents ${a_{i_j}}^{rN'}$. This implies $\abs{N} =
\abs{rN' + N_j} \leq rN_0 + p$. 

For $d > 1$ write $w$ in the form \textup{(ii)} or \textup{(ii$'$)}. Then
the terminal subword has depth $d-1$ in $u_j$ for some $j$. By the
induction hypothesis $u_j$ represents ${c_{\sigma(i_j)}}^{N'}$ where 
\begin{equation} \label{depthbound2} 
N_0 r^{d-2} - p(r^{d-3} + \cdots + 1) \ \leq \ \abs{N'} \ \leq \ N_0
r^{d-1} + p(r^{d-2} + \cdots + 1).
\end{equation} 
Then $(s_{i_j} u_j s_{i_j}^{-1})$ represents ${a_{i_j}}^{rN'}$ and $rN' - p
\leq \abs{N} \leq rN' + p$. These bounds and \eqref{depthbound2} together
imply \eqref{depthbound}. 
\end{proof}

\begin{proposition}[Snowflake word length] \label{sflength} 
Given $r$ and $P$ there are positive constants $C_0, C_1$ with the
following property. If $c$ is the diagonal element of one of the vertex
groups and $w$ is a snowflake word representing $c^N$ then $C_0
\abs{w}^{\alpha} \leq \abs{N} \leq C_1 \abs{w}^{\alpha}$, where $\alpha =
\log_{\lambda}(r)$ and $\lambda$ is the Perron--Frobenius eigenvalue of
$P$. 
\end{proposition}

\begin{proof}
If $w$ is non-trivial and has the form \textup{(i)} or \textup{(i$'$)}
then $1 \leq \abs{N} \leq N_0$ and $\abs{N} \leq \abs{w} \leq
r\abs{N}$. Then $\abs{w}^{\alpha} \leq (rN_0)^{\alpha}$, which implies 
\[ (rN_0)^{-\alpha} \abs{w}^{\alpha} \ \leq \ \abs{N} \ \leq
\ \abs{w}^{\alpha}.\]
Thus we need to arrange that $C_0 \leq (rN_0)^{-\alpha}$ and $C_1 \geq
1$ to cover this case. 

Next assume that $w$ is of type \textup{(ii)} or \textup{(ii$'$)}, which
implies that the depth of every terminal subword is at least
one. Equivalently, $w$ contains the letters $s_j$, $s_j^{-1}$ for some
$j$. Let $s(w)$ be the number of letters $s_j$ or $s_j^{-1}$ in $w$ (for
all indices $j$). Note that a subword of $w$ containing no such letters
has length at most $rN_0$. Since $s(w) \not= 0$, this implies
\begin{equation}\label{swbound} 
s(w) \ \leq \ \abs{w} \ \leq \ 2(rN_0 +1) s(w). 
\end{equation}
Hence $s(w)$ gives an approximate measure of the length of $w$. It can
be computed explicitly, by following the evolution of the curve shortening
process, which in turn is governed by the matrix $P$. 
Note that matched $s$, $s^{-1}$ pairs enclose snowflake subwords
representing powers of $c_j$ for various $j$. These subwords will be
called \emph{$c_j$--subwords}. 

We claim that if $w$ represents a power of $c_k$, and every terminal
subword has depth $i$ or greater, then the number of $c_j$--subwords of
depth $i$ is given by the $kj$--entry of $P^i$, denoted $p^{(i)}_{kj}$.

If $i=1$ then the claim is evident from expressions \textup{(ii)} and
\textup{(ii$'$)}, since the entry $p_{kj}$ of $P$ gives the number of
directed edges from vertex $v_k$ to vertex $v_j$ (and hence the number of
occurrences of $j$ among the indices $\sigma(i_1), \ldots, \sigma(i_m)$). 
Similarly, for $i>1$, each $c_{\ell}$--subword of depth $i-1$ contains
$p_{\ell j}$  $c_j$--subwords of depth $i$, by \textup{(ii)} and 
\textup{(ii$'$)}. The claim now follows by induction on $i$: 
summing over all snowflake subwords of depth $i-1$ and applying the
induction hypothesis, the total number of $c_j$--subwords of depth $i$ is 
$\sum_{\ell} p^{(i-1)}_{k \ell} p_{\ell j} = p^{(i)}_{k j}$. 

Let $x_1, \ldots, x_R$ be the standard basis vectors of $\R^R$. Also 
let $\aabs{ \, \cdot \, }_1$ denote the $\ell_1$ norm on $\R^R$:
$\aabs{v}_1$ is the sum of the absolute values of the entries of the
vector $v$. Let $P^T$ be the transpose of $P$. 

The $kj$--entry of $P^i$ is equal to the $j$--entry of the
column vector $(P^T)^i(x_k)$. Suppose for the moment that every terminal
subword of $w$ has depth $d$. Then for $i \leq d$, the total number of
$s,s^{-1}$ pairs enclosing snowflake subwords of depth $i$ is given by
$\aabs{(P^T)^i(x_k)}_1$. Hence we have 
\[ s(w) \ = \ 2\left(\aabs{P^T(x_k)}_1 + \aabs{(P^T)^2(x_k)}_1 + \cdots +
\aabs{(P^T)^d(x_k)}_1\right).\] 

If we let $d_0$ and $d_1$ denote the smallest and largest depths of
terminal subwords of $w$ then we obtain
\[ 2\sum_{i=1}^{d_0} \aabs{(P^T)^i(x_k)}_1 
\ \leq \ s(w) \ \leq \ 2\sum_{i=1}^{d_1} \aabs{(P^T)^i(x_k)}_1 .\]
Applying \fullref{pfprop} with the norm $\aabs{ \, \cdot \ }_1$
we have
\[ 2A_0 \sum_{i=1}^{d_0} \lambda^i \ \leq \ s(w) \ \leq \ 2A_1 \sum_{i=1}^{d_1}
\lambda^i \ = \ \frac{2A_1 \lambda}{\lambda -1} (\lambda^{d_1} - 1) \]
which implies 
\[ 2A_0 \lambda^{d_0} \ \leq \ s(w) \ \leq \ \frac{2A_1 \lambda}{\lambda
-1} \lambda^{d_1}.\]
Hence by \eqref{swbound} we have 
\begin{equation}
(2A_0) \lambda^{d_0} \ \leq \ \abs{w} \ \leq \ \left(\frac{4(rN_0 + 1)
A_1 \lambda }{\lambda -1}\right) \lambda^{d_1}. 
\end{equation}
We complete the proof by applying \fullref{sfdepth} separately for the
upper and lower bounds. Using $d=d_1$ we obtain 
\[ \abs{N} \ \geq \ B_0 r^{d_1} \ = \ B_0 (\lambda^{d_1})^{\log_{\lambda}(r)}
\ \geq \ B_0 \left(\frac{4(rN_0 + 1) A_1 \lambda }{\lambda
-1}\right)^{-\log_{\lambda}(r)} \abs{w}^{\log_{\lambda}(r)}. \]
Now choose $C_0>0$ satisfying $C_0 \leq B_0 \left(\frac{4(rN_0 + 1) A_1
\lambda }{\lambda -1}\right)^{-\alpha}$ and $C_0 \leq
(rN_0)^{-\alpha}$ to obtain the desired lower bound. 

Applying \fullref{sfdepth} with $d = d_0$ gives
\[ \abs{N} \ \leq \ B_1 r^{d_0} \ = \ B_1
(\lambda^{d_0})^{\log_{\lambda}(r)} \ \leq \ B_1
(2A_0)^{-\log_{\lambda}(r)} \abs{w}^{\log_{\lambda}(r)}\] 
so choose $C_1$ with  $C_1 \geq B_1 (2A_0)^{-\alpha}$ and $C_1 \geq 1$. 
\end{proof}

\section{Proof of Theorem A} \label{thm1sect}

Throughout this section $G_{r,P}$ is fixed, with $r = p/q$ greater than all 
the row sums of $P$,  and $\alpha =
\log_{\lambda}(r)$, where $\lambda$ is the Perron--Frobenius eigenvalue
of $P$. Unless otherwise stated, all words use the generating set $\{a_1,
\ldots, a_n, s_1, \ldots, s_n\}$ for $G_{r.P}$. 

\subsection*{The lower bound} 
To establish the lower bound $\delta(x) \succeq x^{2\alpha}$ we will show
that $\delta (n_i) \geq ({C_0}^2 4^{-\alpha}) \, {n_i}^{2\alpha}$ for
certain integers $n_i$ tending to infinity. This is sufficient by
\fullref{sparse}, provided the sequence $(n_i)$ grows at most 
exponentially. 

Note also that to establish a single inequality $\delta(n) \geq A$, it
is enough to exhibit an embedded disk in $\widetilde{X}_{r,P}$ with
boundary length $n$ and area $A$ or greater, by \fullref{embedded}. Here
we are using the facts that $X_{r,P}$ is aspherical and $2$--dimensional. 

Choose a vertex group $V_m$ in $G_{r,P}$ with $m \geq 2$ and let $c$ be
its diagonal element. There must be at least one vertex group of this
type, for otherwise $P$ would be a permutation matrix with
Perron--Frobenius eigenvalue $1$. For each $i$ choose positive and
negative snowflake words $w_i^+$ and $w_i^-$ representing $c^i$. Then
define $w_i = w_i^+ (w_i^-)^{-1}$ and $n_i = \abs{w_i}$. Note that $C_0
\, 2^{-\alpha} \abs{w_i}^{\alpha} \leq i \leq C_1\, 2^{-\alpha}
\abs{w_i}^{\alpha}$ by \fullref{sflength}. It follows that the
sequence $(n_i)$ tends to infinity, and that it is exponentially bounded: 
\begin{equation*}
\frac{n_{i+1}}{n_i}
\ \leq \ \Biggl(\frac{(i+1)C_1}{iC_0}\Biggr)^{1/\alpha} \
\leq \ \Biggl(\frac{2C_1}{C_0}\Biggr)^{1/\alpha} 
\end{equation*}
for $i \geq 1$. 

Next we find embedded disks $\Delta_i$ in $\widetilde{X}_{r,P}$ with
boundary words $w_i$ and estimate their areas. Each $\Delta_i$ is made of
two disks $\Delta_i^+$ and $\Delta_i^-$ with boundary words $w_i^+ c^{-i}$
and $c^i(w_i^-)^{-1}$ respectively, joined along the boundary arcs labeled 
$c^{-i}, c^i$. After joining, the arc labeled $c^i$ will be called the
\emph{diameter} of $\Delta_i$. 

\begin{figure}[ht!]
\begin{center}
\labellist
\hair 2pt
\pinlabel {\scriptsize $c^N$} [b] at 217 223
\endlabellist
\includegraphics[width=3.5in]{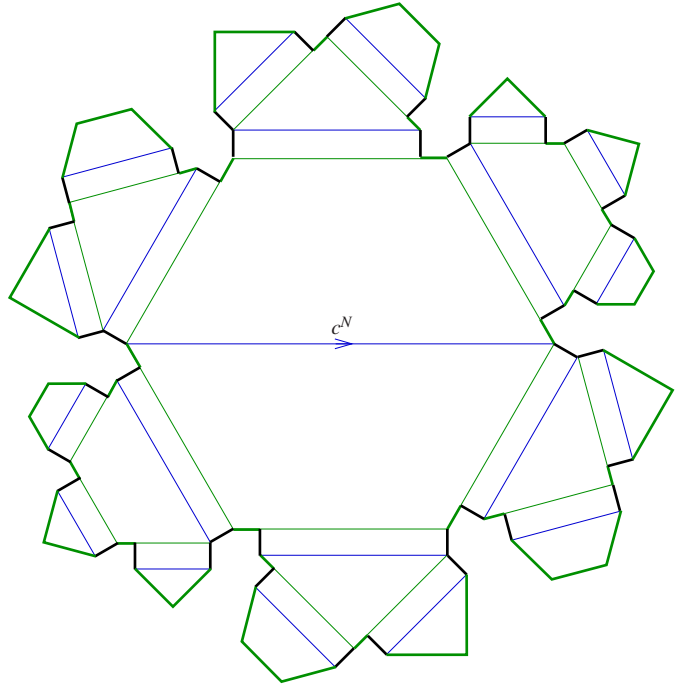}
\end{center}
\caption{A snowflake disk based on the matrix $P =
\bigl(\, {}^1_2 \, {}^1_1\, \bigr)$. The upper and lower halves of the
boundary curve are positive and negative snowflake paths representing
$c^N$.} \label{fig:sfdisk} 
\end{figure}

The disk $\Delta_i^{\pm}$ is a union of embedded disks in vertex spaces
$\widetilde{X}_{m_i}$ and pieces of strips joining them. Consider the
curve shortening process that transforms $c^i$ into $w_i^{\pm}$. To 
build $\Delta_i^{\pm}$ simply fill the central region shown in
\fullref{fig:curveshortening} with the embedded disk from
\fullref{fig:relation}(b). Then fill each strip with either $\lfloor i/p 
\rfloor$ or $\lceil i/p \rceil$ copies of the $2$--cell with the
appropriate boundary word $s_j {c_{\sigma(j)}}^q s_j^{-1} a_j^{-p}$, and
repeat the procedure. The resulting disk is a union of embedded disks in
$\widetilde{X}_{r,P}$ joined along boundary arcs, with no folding along
these arcs. Since each strip separates $\widetilde{X}_{r,P}$, one can see
inductively (on the number of strips crossed by $\Delta_i^{\pm}$) that
$\Delta_i^{\pm}$ is embedded. For the same reason, it suffices to note
that no folding occurs when $\Delta_i^+$ and $\Delta_i^-$ are joined
together to conclude that $\Delta_i$ is embedded. \fullref{fig:sfdisk}
shows an example of a disk $\Delta_i$ with boundary word $w_i$. 

To estimate the area of $\Delta_i$ consider the central region in
$\Delta_i^+$ adjacent to $\Delta_i^-$. By
\fullref{Vmdisks} this subdisk of $\Delta_i$ has area $(m-1) i^2 \geq
i^2$. Then since $i \geq C_0 \, 2^{-\alpha} 
{n_i}^{\alpha}$ (as observed above) we conclude that 
\begin{equation}\label{sfdisklowerbound}
\area(\Delta_i) \geq ({C_0}^2 4^{-\alpha}) {n_i}^{2\alpha} 
\end{equation}
and therefore $\delta(n_i)
\geq ({C_0}^2 4^{-\alpha}) {n_i}^{2\alpha}$.

\subsection*{The upper bound} 
Suppose a word $w$ represents an element of a vertex group $V_m$. The
graph of groups structure of $G_{r,P}$ yields a decomposition of $w$ as
$w_1 \cdots w_k$ where each $w_i$ is either an element of $V_m$, 
or begins with $s_j^{\pm}$ and ends with $s_j^{\mp}$ for some $j$. These
latter cases occur when the path described by $w$ leaves the vertex space
$\widetilde{X}_m$ and then returns again over a strip in
$\widetilde{X}_{r,P}$. 

Recall that a strip in $\widetilde{X}_{r,P}$ has
sides labeled $a_i$ and $c_{\sigma(i)}$. The next lemma shows that a
geodesic (in the generators $\{a_1, \ldots, a_n, s_1, \ldots, s_n\}$) can
only enter a strip from (and return to) the $a_i$--side. 

\begin{lemma}\label{geostrip} 
Let $w$ be a geodesic in $G_{r,P}$ representing an element of a vertex
group $V_m$. Then $w$ is a product of subwords $w_1 \cdots w_k$ where
each $w_i$ is a power of a generator $a_j$, or begins with $s_j$ and ends
with $s_j^{-1}$ (for some $j$) and represents a power of $a_j$.  
\end{lemma}

\begin{proof} 
Let $w' \subset w$ be an innermost word that begins with $s_{\ell}^{-1}$
and ends with $s_{\ell}$ (for some $\ell$) and whose corresponding path in
$\widetilde{X}_{r,P}$ has endpoints in the same vertex space
$\widetilde{X}_{v_{\sigma(\ell)}}$. Thus $w' = s_{\ell}^{-1} u s_{\ell}$
crosses a strip from the $c_{\sigma(\ell)}$--side, and the subword $u$
only crosses strips from (and returns to) $a_i$--sides. That is, $u$ can
be written as $u_1 \cdots u_k$ where each $u_i$ is a power of a generator
$a_j$, or begins with $s_j$ and ends with $s_j^{-1}$ and represents a
power of $a_j$. 

Note that $u$ has both endpoints on an $a_{\ell}$--line in the vertex space
$\widetilde{X}_{v_{\rho(\ell)}}$ across a strip from
$\widetilde{X}_{v_{\sigma(\ell)}}$. Hence $u$ represents $a_{\ell}^N$
for some $N$. Let $u'$ be the word in the standard generators of
$G_{v_{\rho(\ell)}} \cong V_m$ obtained by replacing each $u_i$ by the
appropriate power of $a_j$ that it represents. Consider the word $u'
a_{\ell}^{-N}$ which represents the trivial element $c^0$ in $V_m$. Since
$u'$ does not involve $c$, \fullref{shuffle} implies that every
$a_{j}$--exponent of $u' a_{\ell}^{-N}$ is zero. Hence $u'$ has
$a_{\ell}$--exponent $N$ and $a_j$--exponent zero for every $j \not= \ell$. 

If any of the subwords $u_i$ of $u$ represent a power of $a_j$ with $j
\not= \ell$, then by \fullref{shuffle} one could rearrange the subwords
(preserving the property that $u$ represents $a_{\ell}^N$) so that those 
representing powers of $a_j$ are adjacent. Then these adjacent subwords
cancel in $V_m$ and can be deleted, shortening $w$. Therefore every $u_i$
represents a power of $a_{\ell}$. 

If none of the subwords $u_i$ begins with $s_{\ell}$ and ends with
$s_{\ell}^{-1}$ then $u = a_{\ell}^N$, but then $w'$ could be
replaced by a word $a_{i_1}^{N/r} \cdots a_{i_m}^{N/r}$ representing
${c_{\sigma(\ell)}}^{N/r}$. The new word is shorter than $w$ because 
of the hypothesis that $m<r$, and therefore some $u_i$ must have the form
$s_{\ell} v s_{\ell}^{-1}$ after all. Now rearrange the subwords so that
$s_{\ell} v s_{\ell}^{-1}$ occurs last. Again $w$ can be shortened by
replacing $u$ with this rearranged word and then cancelling
$s_{\ell}^{-1} s_{\ell}$ at the end. 
\end{proof}

\begin{proposition}\label{sfgeo}
Let $c$ be the diagonal element of one of the vertex groups in
$G_{r,P}$. Then for every $N$ there is a snowflake word $w_{\text{sf}}$
and a geodesic $w_{\text{geo}}$, both representing $c^N$, with
$\abs{w_{\text{sf}}} \leq r N_0 \abs{w_{\text{geo}}}$. 
\end{proposition}

\begin{proof} The proof is by induction on $\abs{N}$. 
Let $w$ be a geodesic representing $c^N$. We shall apply
\fullref{shuffle} inductively to rearrange and modify $w$ into two words,
a geodesic $w_{\text{geo}}$ and a positive snowflake word
$w_{\text{sf}}$. The two constructions are identical except at the base 
of the induction, which involves only certain segments of length at most
$rN_0$. 

Let $a_{i_1}, \ldots, a_{i_m}$ be the standard generators (in order) of
the vertex group $V_m$ containing $c$.  If $\abs{N} \leq N_0$ then
define $w_{\text{geo}} = w$ and $w_{\text{sf}}=a_{i_1}^N \cdots
a_{i_m}^N$. The desired conclusion holds in this case since $r>m$. 

Suppose next that $\abs{N} > N_0$. By \fullref{geostrip} we can write
$w$ as $w_1 \cdots w_k$ where each subword has the form
$a_{{j}}^{N_{j}}$ or $s_{{j}} u_{j} s_{{j}}^{-1}$. In the latter case
$s_{{j}} u_{j} s_{{j}}^{-1}$ represents a power of $a_{{j}}$. 

By \fullref{shuffle} we can permute the subwords $w_{\ell}$ of
$w$ to arrange that those representing powers of $a_{i_1}$
come first, those representing powers of $a_{i_2}$ occur next, and so
on. The resulting word is still a geodesic representing $c^N$. Note that
two subwords cannot both be of the form $s_{i_j}
u_{j} s_{i_j}^{-1}$ since they could be made adjacent, and then a
cancellation of $s_{i_j}^{-1} s_{i_j}$ would be possible. 
Hence we can arrange
for $w$ to have the form 
\begin{equation} \label{rearranged} 
 w \ = \ (s_{i_1} u_1 s_{i_1}^{-1}) (a_{i_1}^{N_1})(s_{i_2} u_2
s_{i_2}^{-1})(a_{i_2}^{N_2}) \, \cdots \, (s_{i_m} u_m
s_{i_m}^{-1})(a_{i_m}^{N_m})
\end{equation} 
where each $s_{i_j} u_j s_{i_j}^{-1}$ represents a power
of $a_{i_j}$. 
Next observe that $\abs{N_{j}}<p$ for all $j$,
since otherwise a subword of the form $s_{i_j}^{-1}
a_{i_j}^{\pm p}$ could be replaced by 
a word of the form $a_{\ell_1}^{\pm q} \cdots a_{\ell_{m'}}^{\pm q}
s_{i_j}^{-1}$ (that is, ${c_{\sigma(i_j)}}^{\pm q} s_{i_j}^{-1}$
expressed in the standard generators). Here $m'$ is a row sum of $P$ and
so $r > m'$, making the new word shorter than $w$. 

Recall that $u_j$ represents a power of $c_{\sigma(j)}$. 
By \fullref{shuffle} the power of $a_{i_j}$ represented by $s_{i_j} u_j
s_{i_j}^{-1}$ is $N - N_j$, and so $u_j$ represents ${c_{\sigma(j)}}^{(N
- N_j)/r}$. Recall that $N_0 > p$, hence $\abs{N} > p >
\abs{N_j}$. Then since $r > 2$ it follows that $\abs{(N - N_j)/r} <
\abs{N}$.

By induction ${c_{\sigma(j)}}^{(N - N_j)/r}$ is represented by a geodesic
$(u_j)_{\text{geo}}$ and a positive snowflake word $(u_j)_{\text{sf}}$
satisfying the conclusion of the lemma. Define $w_{\text{geo}}$ and
$w_{\text{sf}}$ by replacing each subword $u_j$ in \eqref{rearranged} by 
$(u_j)_{\text{geo}}$ or $(u_j)_{\text{sf}}$ accordingly. Then the desired
conclusion also holds for $w_{\text{geo}}$ and $w_{\text{sf}}$, since
they agree except in the subwords $(u_j)_{\text{geo}}$ and
$(u_j)_{\text{sf}}$. 
\end{proof}

\begin{corollary}[Edge group distortion]\label{distortion} 
Given $r$ and $P$ there is a positive constant ${D}$ with the following
property. If $c$ is a diagonal element and $w$ is a word representing
$c^N$ then $\abs{N} \leq {D} \abs{w}^{\alpha}$. 
\end{corollary}

\begin{proof}
It suffices to consider the case when $w$ is a geodesic. Apply 
\fullref{sfgeo} to obtain the geodesic $w_{\text{geo}}$ and
snowflake word $w_{\text{sf}}$ representing $c^N$ with
$\abs{w_{\text{sf}}} \leq r N_0 \abs{w_{\text{geo}}}$. Then
\fullref{sflength} implies $\abs{N} \leq C_1 \abs{w_{\text{sf}}}^{\alpha} 
\leq C_1 \left( r N_0\right)^{\alpha} \abs{w_{\text{geo}}}^{\alpha}$. 
\end{proof}

The statement and proof of the next proposition are similar to 
those of Brady--Bridson \cite[Proposition 3.2]{BrBr}. The case $N = 0$
establishes the upper bound of \fullref{mainthm1}. 

\begin{proposition}[Area bound]\label{toplevel2} 
Given $r$ and $P$ there is a positive constant ${E}$ with the following
property. If $w$ is a word in $G_{r, P}$ representing $x^N$ for some $N$,
where $x$ is either a generator $a_i$ or the diagonal element of one of
the vertex groups, then $\area(w x^{-N}) \leq {E} \abs{w}^{2\alpha}$. 
\end{proposition}

\proof
We argue by induction on $\abs{w}$. We shall prove the statement with
${E} = (3/2)r^2 {D}^2$ (${D}$ given by \fullref{distortion}). 
Let $c$ denote the diagonal element
of the vertex group $V_m$ containing $x$. 

Write $w$ as $w_1 \cdots w_k$ where each $w_i$ has the form
$a_{j_i}^{N_i}$ or is a word beginning in $s_{j_i}^{\pm 1}$ and ending in
$s_{j_i}^{\mp 1}$. In the latter cases $w_i$ represents an element of the
form $c^{N_i}$ or $a_{j_i}^{N_i}$. Let $I_c$ and $I_a$ be the sets of
indices for which these two cases occur, and let $w'$ be the word
obtained from $w$ by replacing each subword $w_i$ of this type with the
appropriate word $c^{N_i}$ or $a_{j_i}^{N_i}$. Then $w'$ is a word in the
standard generators of $V_m$ (and the diagonal element) representing
$x^N$, of length $\sum_i N_i$. 

By \fullref{Vmarea} we have $\area(w' x^{-N}) \leq 3 \sum_{i<j} N_i
N_j$. To estimate each $N_i$ we use \fullref{distortion} as
follows. If $i \in I_c$ then $w_i$ represents $c^{N_i}$ and
\fullref{distortion} gives $N_i \leq {D} \abs{w_i}^{\alpha}$. If $i \in
I_a$ then $w_i = s_{j_i} u_i s_{j_i}^{-1}$ for some $u_i$ representing
${c_{\sigma(j_i)}}^{N_i / r}$ (because $w_i$ represents
$a_{j_i}^{N_i}$). Then by \fullref{distortion} we have $N_i/r
\leq {D}(\abs{w_i}-2)^{\alpha} \leq {D}\abs{w_i}^{\alpha}$, so $N_i \leq
r{D}\abs{w_i}^{\alpha}$. Finally if $i \not\in (I_c \cup I_a)$ then $N_i =
\abs{w_i} \leq \abs{w_i}^{\alpha}$. Putting these observations together
we have 
\begin{equation}\label{basearea}
\area(w' x^{-N}) \ \leq \ 3 r^2 {D}^2 \sum_{i<j} \abs{w_i}^{\alpha}
\abs{w_j}^{\alpha}. 
\end{equation}

Next we use the induction hypothesis and \fullref{distortion} to
bound $\area(w w'^{-1})$. First note that $\area(w w'^{-1}) \ \leq \ 
\sum_{i\in I_c} \area(w_i {c}^{-N_i}) + \sum_{i \in I_a} \area(w_i
a_{j_i}^{-N_i})$. 

If $i \in I_c$ then $w_i = s_{j_i}^{-1} u_i s_{j_i}$ where $u_i$
represents ${a_{j_i}}^{rN_i}$. Applying the induction hypothesis to $u_i$
we have $\area(u_i {a_{j_i}}^{-rN_i}) \leq (3/2) r^2 {D}^2 (\abs{w_i} -
2)^{2\alpha}$. The strip $s_{j_i}^{-1} {a_{j_i}}^{rN_i} s_{j_i} c^{-N_i}$
has area $N_i/q \leq ({D}/q) \abs{w_i}^{\alpha} \leq {D}\abs{w_i}^{\alpha}$,
by \fullref{distortion}. Thus 
\begin{equation}\label{wicarea}
\begin{split}
\area(w_i c^{-N_i}) \ &\leq \ (3/2) r^2 {D}^2 (\abs{w_i} -2)^{2\alpha} +
{D}\abs{w_i}^{\alpha} \\
&\leq \ (3/2) r^2 {D}^2 ( (\abs{w_i} - 2)^{2\alpha} + \abs{w_i}^{\alpha})
\\ 
&\leq \ (3/2) r^2 {D}^2 \abs{w_i}^{2\alpha}. 
\end{split}
\end{equation}
The last inequality above uses the fact that for numbers $x \geq 0$ one
has $(x+2)^{2\alpha} \geq x^{\alpha}(x+2)^{\alpha} +
2^{\alpha}(x+2)^{\alpha} \geq x^{2\alpha} + (x+2)^{\alpha}$. 

If $i \in I_a$ then $w_i = s_{j_i} u_i s_{j_i}^{-1}$ where $u_i$
represents ${c_{\sigma(j_i)}}^{N_i / r}$.  Applying the induction
hypothesis to $u_i$ we have $\area(u_i {c_{j_i}}^{-N_i/r}) \leq (3/2) r^2
{D}^2 (\abs{w_i} - 2)^{2\alpha}$. The strip $s_{j_i} {c_{j_i}}^{N_i/r}
s_{j_i}^{-1} a_{j_i}^{-N_i}$ has area $(N_i/r)/q \leq ({D}/q)
(\abs{w_i}-2)^{\alpha} \leq {D}(\abs{w_i} - 2)^{\alpha}$, by
\fullref{distortion}. Therefore 
\begin{equation}\label{wiaarea}
\begin{split}
\area(w_i a_{j_i}^{-N_i}) \ &\leq \ (3/2) r^2 {D}^2 (\abs{w_i} -2)^{2\alpha} +
{D}(\abs{w_i} - 2)^{\alpha} \\
&\leq \ (3/2) r^2 {D}^2 ( (\abs{w_i} - 2)^{2\alpha} + (\abs{w_i} - 2)^{\alpha})
\\ 
&\leq \ (3/2) r^2 {D}^2 \abs{w_i}^{2\alpha}. 
\end{split}
\end{equation}
Combining \eqref{wicarea} and \eqref{wiaarea} we then have  
\begin{equation}\label{indhyparea}
\area(w w'^{-1}) \ \leq \ \sum_{i \in I_c \cup I_a} (3/2) r^2 {D}^2
\abs{w_i}^{2\alpha} \ \leq \ \sum_i (3/2) r^2 {D}^2 \abs{w_i}^{2\alpha}. 
\end{equation}
Finally, adding \eqref{basearea} and \eqref{indhyparea} together gives
the desired result: 
\begin{multline*}
\area(w x^{-N}) \ \leq \ (3/2) r^2 {D}^2 \Bigl(\sum_i
\abs{w_i}^{\alpha}\Bigr)^2 \\ \leq \ (3/2) r^2 {D}^2 \Bigl(\sum_i
\abs{w_i}\Bigr)^{2\alpha} \ = \ (3/2) r^2 {D}^2
\abs{w}^{2\alpha}. \quad \qed
\end{multline*}
%
%

\section{Suspension and snowflake balls} \label{s:suspension} 

Throughout this section $P$ denotes a non-negative $R\times R$ integer
matrix with  Perron--Frobenius eigenvalue $\lambda$, and $r$ 
is an integer which is strictly greater than the largest row sum of $P$. 
In this section, we give an explicit description of the suspended 
snowflake groups $\Sigma G_{r,P}$ and the $3$--dimensional $K(\Sigma
G_{r,P}, 1)$ spaces $X^{3}_{r,P}$. Then we describe \emph{snowflake
balls} $B^{3}_{i}$ which embed in the universal cover of $X^{3}_{r,P}$
and estimate their boundary areas. We show how to iterate this 
suspension procedure to obtain groups $\Sigma^{k} G_{r,P}$ and
$(k+2)$--dimensional spaces $X^{k+2}_{r,P}$. Lastly we define
higher-dimensional snowflake balls and estimate their boundary volumes. 

\begin{remark}
In order to realize the exponents $(k + 1)/k$ (the endpoints of the
intervals in \fullref{fig:spectra}, which are omitted otherwise) we
add the free abelian group ${\Z}^{2}$ to the class of snowflake
groups $G_{r,P}$. We endow ${\Z}^{2}$ with snowflake structure 
as follows
\[ {\Z}^{2}\; = \; \langle \, a_{1}, a_{2}, c \; | \; 
  a_{1}a_{2}=c, \ c =a_{2}a_{1} \, \rangle \] 
and use the corresponding presentation $2$--complex $X$ in place of
$X_{r,P}$. There is no matrix $P$ associated to the group ${\Z}^{2}$, and
so the only condition that we impose on the integer $r$ is 
that $r \geq 2$. Since there are no stable letters $s_{i}$, we define the
snowflake words to be the commutators $w_{i} \; = \; [a_{1}^{r^{i}},
a_{2}^{r^{i}}]$ and define the snowflake disks $B^{2}_{i} \; = \;
\Delta_{r^{i}} $ to be the unique embedded disks in $X$ with boundary
$w_i$. 

In the discussions that follow, whenever we talk about 
snowflake groups $G_{r,P}$, we shall always include 
${\Z}^{2}$, and whenever we use the complexes 
$X_{r,P}$ we shall always include the presentation 
$2$--complex $X$ for ${\Z}^{2}$ described above.  
\end{remark}

\subsection*{The groups $\Sigma G_{r,P}$} 
Let $\phi \co G_{r,P} \to G_{r,P}$ be the monomorphism 
which takes each $a_i$ to $a_i^r$ and each $s_i$ to itself. 
The group $\Sigma  G_{r,P}$ is defined to be the associated multiple  
HNN extension with stable letters $u_{1}$ and $v_{1}$: 
\[ \Sigma G_{r,P} \ = \ \langle \, G_{r,P}, u_{1}, v_{1} \mid 
u_{1}gu_{1}^{-1} = \phi (g), \,  v_{1}gv_{1}^{-1}=\phi (g) \;
(g \in G_{r,P}) \, 
\rangle. \]

\subsection*{The spaces $X^3_{r,P}$} These spaces will have fundamental
group $\Sigma G_{r,P}$. Recall that $X_{r,P}$ is a 
$2$--dimensional $K(G_{r,P}, 1)$ space. There is a cellular map
$\Phi \co X_{r,P} \to X_{r,P}$ which induces the map $\phi$ on the
fundamental group. It maps the $1$--cells labeled $s_{i}$ 
homeomorphically to themselves, maps the $1$--cells labeled 
$a_{i}$ to themselves by degree $r$ maps, and maps each 
$2$--cell in the obvious manner;  the image 
of each triangular $2$--cell has combinatorial area  $r^{2}$, and the 
image of the remaining $2$--cells (which have an $s_{i}$ edge  in their 
boundaries) have combinatorial area $r$.   The $3$--complex $X^3_{r,P}$ 
with fundamental group  $\Sigma G_{r,P}$ is obtained by taking two copies 
of the mapping torus of the map $\Phi$ and identifying them along a copy of 
$X_{r,P}$. From this perspective it is easy to see that 
$X^3_{r,P}$ is aspherical; each mapping torus 
is aspherical since $X_{r,P}$ is an aspherical $2$--complex, 
and since $\Phi$ induces the monomorphism $\phi$ in $\pi_{1}$. 
We give more details of the cell structure of $X^3_{r,P}$ below. 

Start with the $2$--complex $X_{r,P}$ and form two copies of $X_{r,P}
\times [0,1]$. Each copy is given the product cell structure, in which
each $k$--cell of $X_{r,P}$ gives rise to a $(k+1)$--cell in $X_{r,P} \times
(0,1)$. The ``bottom'' side $X_{r,P} \times \{0\}$ keeps its original
cell structure and the ``top'' $X_{r,P} \times \{1\}$ is subdivided by
pulling back under $\Phi$ the cell structure of $\Phi(X_{r,P})$. That
is, each triangular $2$--cell in a vertex space of $X_{r,P}$ is subdivided
into $r^2$ triangles, and each edge space $2$--cell (bearing the boundary
label $s_j c_{\sigma(j)} s_j^{-1} a_j^r$) is subdivided into $r$ copies
of the same cell. 

\begin{figure}[ht!]
	\begin{center}
\labellist
\hair 2pt
\pinlabel* {\scriptsize $a$} [Bl] at 82 192
\pinlabel* {\scriptsize $a$} [Bl] at 298 192
\pinlabel* {\scriptsize $a$} [Bl] at 150 124
\pinlabel* {\scriptsize $a$} [Bl] at 109 75
\pinlabel* {\scriptsize $b$} [Bl] at 145 192
\pinlabel* {\scriptsize $b$} [Bl] at 362 192
\pinlabel* {\scriptsize $b$} [Bl] at 214 124
\pinlabel* {\scriptsize $b$} [Bl] at 293 70
\pinlabel* {\scriptsize $c$} [Bl] at 134 240
\pinlabel* {\scriptsize $c$} [Bl] at 350 240
\pinlabel* {\scriptsize $c$} [Bl] at 201 172
\pinlabel* {\scriptsize $v_{1}$} [Bl] at 5 183
\pinlabel* {\scriptsize $v_{1}$} [Bl] at 482 183
\pinlabel* {\scriptsize $v_{1}$} [Bl] at 185 48
\endlabellist
		\includegraphics[width=2.5in]{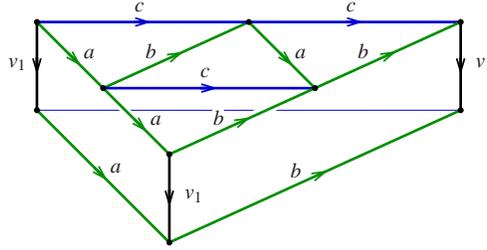}
	\end{center}
	\caption{A triangular $3$--cell (with $r=2$)}
	\label{fig:prism}
\end{figure}

\begin{figure}[ht!]
	\begin{center}
\labellist
\hair 2pt
\pinlabel {\scriptsize $a_j$} [bl] at 84 83
\pinlabel {\scriptsize $a_j$} [bl] at 183 83
\pinlabel {\scriptsize $a_j$} [bl] at 291 83
\pinlabel {\scriptsize $a_j$} [bl] at 399 83
\pinlabel {\scriptsize $a_j$} [bl] at 130 11
\pinlabel {\scriptsize $a_j$} [bl] at 345 11
\pinlabel {\scriptsize $s_j$} [br] at 47 100
\pinlabel {\scriptsize $s_j$} [br] at 263 100
\pinlabel {\scriptsize $s_j$} [br] at 479 100
\pinlabel {\scriptsize $s_j$} [tl] at 495 35
\pinlabel {\scriptsize $c_{\sigma(j)}$} [b] at 190 139
\pinlabel {\scriptsize $c_{\sigma(j)}$} [b] at 407 139
\pinlabel {\scriptsize $v_{1}$} [r] at 22 40
\pinlabel {\scriptsize $v_{1}$} [r] at 239 40
\pinlabel {\scriptsize $v_{1}$} [r] at 455 40
\pinlabel {\scriptsize $v_{1}$} [l] at 522 94
\endlabellist
		\includegraphics[width=3.0in]{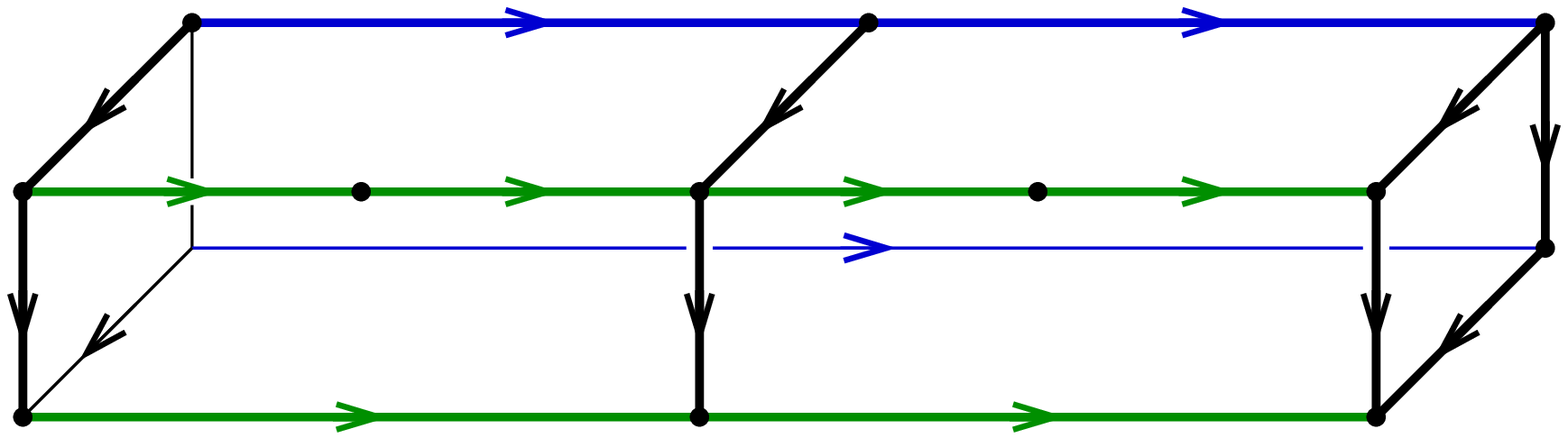}
	\end{center}
	\caption{A rectangular $3$--cell}
	\label{fig:prism2}
\end{figure}

The vertical $1$--cells of the two copies of $X_{r,P} \times [0,1]$ are
labeled $u_{1}$ and $v_{1}$ respectively, oriented from $X_{r,P} \times
\{1\}$ to $X_{r,P} \times \{0\}$. Finally to form $X^3_{r,P}$ one
attaches the bottom of each piece to $X_{r,P}$ by the identity, and the
top by the map $\Phi$. Figures \ref{fig:prism} and \ref{fig:prism2}
illustrate the two types of $3$--cell occurring in $X^3_{r,P}$. 

\subsection*{Snowflake balls}
We define embedded $3$--dimensional balls $B^{3}_j$ in
$\widetilde{X}^3_{r,P}$ in a similar fashion to the snowflake disks
constructed in \fullref{thm1sect}. An essential difference, however,
is that now $r$ is an integer, and the observations of
\fullref{integerpower} apply. That is, snowflake disks
with diameter labeled $c^{r^i}$ are unique, and the corresponding
snowflake words have no ``remainder'' terms. 

As in the proof of \fullref{mainthm1} we let $c$ be the diagonal
element of a vertex group $V_m$ in $G_{r,P} \subset \Sigma G_{r,P}$ where $m
\geq 2$. We let $w_{i}^{+}$ and $w_{i}^{-}$ denote respectively the (unique)
positive and negative snowflake words representing $c^{r^{i}}$. (Note
that the indexing here differs from that in \fullref{thm1sect}, where
these words would be called $w_{r^i}^{\pm}$.) Let $B^{2}_{i}$ be the
snowflake disk bounded by $w_i = w_i^+ (w_i^-)^{-1}$, with diameter
labeled $c^{r^i}$. Note that $B^{2}_{i}$ is the same as the snowflake  
disk $\Delta_{r^{i}}$ of \fullref{thm1sect}. 

For each positive integer $j$, we shall use a {\em stack of thickened van
Kampen disks} to define an embedded $3$--ball $B^{3}_{j}$ in the
universal cover of $X^3_{r,P}$.  Note that the universal cover of
$X^3_{r,P}$ contains infinitely many embedded copies of the universal
cover of $X_{r,P}$; one for each coset of $G_{r,P}$ in $\Sigma
G_{r,P}$. We call two such copies {\em adjacent} if the cosets have
representatives which differ by right multiplication by $u_{1}^{\pm 1}$
or $v_{1}^{\pm 1}$.

The map $\Phi\co X_{r,P} \to X_{r,P}$ lifts to a map
of universal covers which we also denote by $\Phi$. Consider the image 
${\Phi}(B^{2}_{i})$ of the embedded snowflake disk
$B^{2}_i$. This image is again embedded, but its boundary word is
$\phi(w_{i})$. If we apply the curve shortening procedure once to the
subword $\phi(w_{i}^{+})$ we obtain $w_{i+1}^{+}$, which is the
positive snowflake word for $c^{r^{i+1}}$. Similarly, if we apply curve
shortening once to the subword $\phi(w_{i}^{-})$ we obtain the negative
snowflake word for $c^{r^{i+1}}$. Thus ${\Phi}(B^{2}_{i})$
is a sub-diagram of $B^{2}_{i+1}$. The {\em top half} of the ball $B^{3}_{j}$
is defined to be the union of the mapping cylinders of
${\Phi}$ with domain $B^{2}_{i}$ and codomain $B^{2}_{i+1}$
where $i$ ranges from 1 to $j$; the copies of $B^{2}_{i}$ are
identified.  This embeds in the universal cover of $X^3_{r,P}$ as follows.
The disk $B^{2}_{1}$ embeds in some copy of the universal cover of
$X_{r,P}$, $B^{2}_{2}$ embeds in the adjacent copy obtained by right
multiplying by $u_{1}^{-1}$, and the mapping cylinder of $\Phi\co
B^{2}_{1} \to B^{2}_{2}$ embeds in the universal cover of $X^3_{r,P}$ to
interpolate between the images of $B^{2}_{1}$ and $B^{2}_{2}$. Note
that this embedding is possible since the universal covering of $X^3_{r,P}$
can be described as an infinite union of mapping cylinders of
${\Phi}\co \widetilde{X}_{r,P} \to \widetilde{X}_{r,P}$ which
is encoded by the Bass--Serre tree $T$ corresponding to the multiple HNN
description of $\Sigma G_{r,P}$. 

\begin{figure}[ht!]
	\begin{center}
\labellist
\hair 2pt
\pinlabel {\scriptsize $u_{1}$} [r] at 59 100
\pinlabel {\scriptsize $u_{1}$} [r] at 96 120
\pinlabel {\scriptsize $v_{1}$} [r] at 59 77
\pinlabel {\scriptsize $v_{1}$} [r] at 96 56
\pinlabel {\scriptsize $B^{2}_1$} [l] at 181 172
\pinlabel {\scriptsize $B^{2}_1$} [l] at 181 3
\pinlabel {\scriptsize $B^{2}_{j+1}$} [l] at 332 87
\pinlabel {\scriptsize $B^{2}_{j}$} [l] at 279 107
\pinlabel {\scriptsize $B^{2}_{j}$} [l] at 279 67
\pinlabel {\scriptsize $B^{2}_{j-1}$} [l] at 243 128
\pinlabel {\scriptsize $B^{2}_{j-1}$} [l] at 243 47
\endlabellist
		\includegraphics[width=4in]{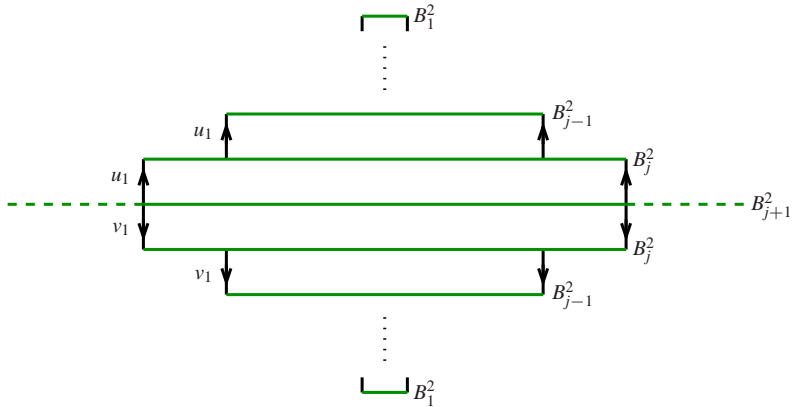}
	\end{center}
	\caption{A schematic diagram of the embedded ball 
	$B^{3}_j$}\label{fig:3ball}
\end{figure}

\begin{figure}[ht!]
	\begin{center}
		\includegraphics[width=3.5in]{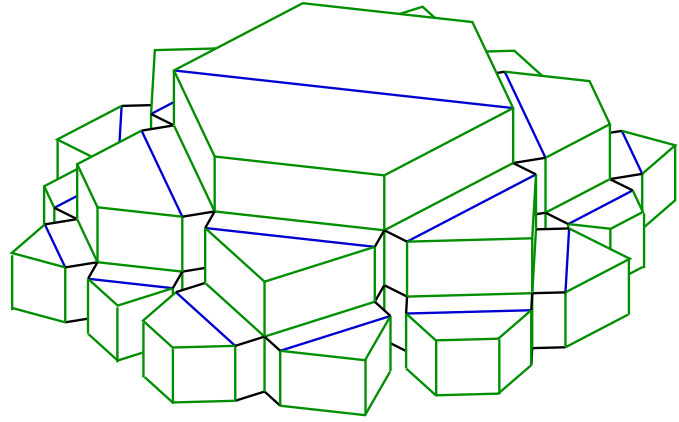}
	\end{center}
	\caption{A few layers of $B^{3}_{j}$}
	\label{fig:sfball}
\end{figure}

We continue to add mapping cylinders of $\Phi\co B^{2}_{i} \to
B^{2}_{i+1}$ for $i = 2, \ldots, j$, as indicated in the top half of the  
schematic diagram in \fullref{fig:3ball}. 
The image of the union of the first  few embedded layers is shown 
in \fullref{fig:sfball}. 
In a similar fashion, we can embed a second copy 
of the union of mapping cylinders of $\Phi\co B^{2}_{i} \to
B^{2}_{i+1}$. However,  this time  we start from the copy of
$B^{2}_{j}$ in the image of the previous union,  and add the mapping
cylinders in descending order (so $i=j, \ldots, 1$) and require that
new copies of the universal cover of $X_{r,P}$ differ by right
multiplication by $v_{1}^{+1}$.  
The image of this family is indicated in the
lower half of the schematic diagram of  
\fullref{fig:3ball}, and the total union is the embedded ball $B^{3}_{j}$. 
It is easy to see that the union embeds, since each mapping cylinder embeds, 
and distinct mapping cylinders correspond to distinct layers in the 
$3$--complex $\widetilde{X}^3_{r,P}$. These layers are distinct, since they 
map to distinct edges of the Bass--Serre tree $T$. Finally, 
there is a $2$--dimensional ``fringe'' at the equator $B^{2}_{j+1}$ level. 
We remove this fringe  by simply replacing the  
two  embeddings of ${\Phi}\co B^{2}_{j} \to B^{2}_{j+1}$ by
embeddings of  
${\Phi}\co B^{2}_{j} \to {\Phi}(B^{2}_{j})$. 

\begin{lemma} \label{spherearea}
Given $r$ and $P$ there is a positive constant $F_0$ such that
$\abs{\partial B^{2}_{j}} \leq \area(\partial B^{3}_{j})$
$\leq
F_0\abs{\partial B^{2}_{j}} $ for every $j$. 
\end{lemma}

\begin{proof}
The ball $B^{3}_{j}$ is a union of $2j$ mapping cylinders. 
See \fullref{fig:3ball} for a schematic representation. 
Its boundary area is twice the area of the {\em upper hemisphere}. 
This latter area is estimated as follows. 
	
For each $1 \leq i \leq j$, there are $\abs{\partial B^{2}_{i}}$ {\em
vertical} (conjugation by $u_{1}$) $2$--cells, which interpolate between
$\partial B^{2}_{i}$ and $\Phi (\partial B^{2}_{i})$. This proves
the first inequality, $\abs{\partial B^{2}_{j}} \leq \area(\partial
B^{3}_{j})$. 

For each $1\leq i \leq j$ there are {\em horizontal} $2$--cells which
interpolate between $\Phi (\partial B^{2}_{i-1})$ and $\partial
B^{2}_{i}$.  In the case $i=1$ there is no loop $\Phi  (\partial
B^{2}_{0})$, and the horizontal $2$--cells just fill the van Kampen
diagram $B^{2}_{1}$.  For any $i$, the horizontal $2$--cell contribution
to the area is bounded above by $\abs{\partial B^{2}_{i}}$. To see this,
note that the horizontal interpolation is a union of pieces of the form
$s_ja_{i_1} \cdots a_{i_m}s_j^{-1}a_{j}^{-r}$ where $\{a_{1}, \ldots, a_{m}\}$
generates a vertex group $V_{m}$, and the stable letter $s_j$ conjugates
the diagonal element of this vertex group to some generator $a_{j}$ of
$G_{r,P}$. The area of this piece is $m$, and its contribution to
$\abs{\partial B^2_i}$ is $m+2$. 

Counting vertical and horizontal $2$--cells for both hemispheres we
obtain 
\[ \area(\partial B^{3}_{j}) \; \leq \; 4\sum_{i=1}^{j}\abs{\partial
B^{2}_{i}} \, .\]
\fullref{sflength} implies that $\abs{w_{i}^{+}} \; \leq \;
C_{0}^{-1/\alpha}r^{i/\alpha}$ and so 
\[ 4\sum_{i=1}^{j}\abs{\partial B^{2}_{i}} 
\; = \; 
8\sum_{i=1}^{j}\abs{w_{i}^{+}} \; \leq \;
8C_{0}^{-1/\alpha}\sum_{i=1}^{j}(r^{1/\alpha})^{i} \, .\] 
The last term is a geometric series, and so is bounded above by 
$F_0'(r^{1/\alpha})^{j}$ for a positive constant $F_0'$ (independent of
$j$). \fullref{sflength} also gives
$C_{1}^{-1/\alpha}r^{j/\alpha} \leq \abs{w_{j}^{+}}$
and so 
\[ \area(\partial B^{3}_{j}) \; \leq \; F_0' r^{j/\alpha} \; \leq \;
\frac{F_0'}{2}C_{1}^{1/ \alpha}\abs{\partial B^{2}_{j}} \, .\]
Now the desired (second) inequality holds by taking $F_0 = (F_0'
/2)C_{1}^{1/\alpha}$. 
\end{proof}

\subsection*{The inductive supension procedure} 
Having discussed $\Sigma G_{r,P}$ we define further suspensions
$\Sigma^{k} G_{r,P}$ having $(k+2)$--dimensional Eilenberg--MacLane
spaces $X^{k+2}_{r,P}$, and $(k+2)$--dimensional snowflake balls
$B^{k+2}_j \subset \widetilde{X}^{k+2}_{r,P}$. We assume that the group
$\Sigma^{k-1}G_{r,P}$, the space $X^{k+1}_{r,P}$, and snowflake balls
$B^{k+1}_j \subset \widetilde{X}^{k+1}_{r,P}$ have already been
constructed. 

First we define the groups $\Sigma^k G_{r,P}$. Let $\phi_k\co \Sigma^{k-1}
G_{r,P} \to \Sigma^{k-1} G_{r,P}$ be the monomorphism which sends $a_i$
to $a_i^r$ and which leaves fixed the stable letters $s_i$, $u_i$, and
$v_i$. We define $\Sigma^{k} G_{r,P}$ to be the multiple ascending HNN
extension with two stable letters $u_k$ and $v_k$, each acting by $\phi_k$: 
\[ \Sigma^{k}G_{r,P} \ = \ \langle \, \Sigma^{k-1} G_{r,P}, u_{k},
v_{k} \mid 
u_{k} g u_{k}^{-1} = \phi_k(g), \, v_{k} g
v_{k}^{-1} = \phi_k(g) \ (g \in \Sigma^{k-1} G_{r,P}) \, \rangle.\]

Next we define the spaces $X^{k+2}_{r,P}$. The homomorphism $\phi_k$ is
induced by a cellular map $\Phi_{k+1}\co X^{k+1}_{r,P} \to
X^{k+1}_{r,P}$. We define $X^{k+2}_{r,P}$ to be the double mapping torus
with monodromy $\Phi_{k+1}$. That is, take two copies of $X^{k+1}_{r,P}
\times [0,1]$, identify the ``bottom'' sides $X^{k+1}_{r,P} \times \{0\}$
to $X^{k+1}_{r,P}$ by the identity, and attach the ``top'' sides
$X^{k+1}_{r,P} \times \{1\}$ to $X^{k+1}_{r,P}$ by the map
$\Phi_{k+1}$. The vertical $1$--cells of the copies of $X^{k+1}_{r,P}
\times [0,1]$ are labeled $u_k$ and $v_k$ respectively, oriented from
$X^{k+1}_{r,P} \times \{1\}$ to $X^{k+1}_{r,P} \times \{0\}$. The
resulting space $X^{k+2}_{r,P}$ is given a cell structure analogous to
that of $X^3_{r,P}$. As before, $X^{k+2}_{r,P}$ is aspherical,
has dimension $k+2$, and has fundamental group $\Sigma^{k} G_{r,P}$. 

Now we define the higher-dimensional snowflake balls. The map
$\Phi_{k+1}$ lifts to a map $\widetilde{X}^{k+1}_{r,P} \to
\widetilde{X}^{k+1}_{r,P}$ which we will also call $\Phi_{k+1}$. We
define $(k+2)$--dimensional balls $B^{k+2}_j$ of diameter $r^j$ for each
$j$ as unions of mapping cylinders (called {\em layers}) of the map
$\Phi_{k+1}$ restricted to 
the $(k+1)$--dimensional balls $B^{k+1}_i$. These mapping cylinders are
assembled as shown in \fullref{fig:3ball}, with $B^{k+1}_i$ in place
of $B^2_i$. More specifically, we assume inductively that $\Phi_{k+1}$
maps $B^{k+1}_i$ into a subcomplex of $B^{k+1}_{i+1}$ for each $i$. Then
the upper hemisphere of $B^{k+2}_j$ is the union of the mapping cylinders
of $\Phi_{k+1} \co B^{k+1}_i \to B^{k+1}_{i+1}$ where $i$ ranges from $1$
to $j-1$, and the mapping cylinder of $\Phi_{k+1} \co B^{k+1}_j \to 
\Phi_{k+1}(B^{k+1}_{j})$. The lower hemisphere is defined similarly, and the
two are identifed along $\Phi_{k+1}(B^{k+1}_{j})$. Note that the
subspaces $B^{k+1}_i - \Phi_{k+1}(B^{k+1}_{i-1})$ of the domains of these
mapping cylinders lie in the boundary of $B^{k+2}_j$. 

Recall that $\Phi_{k+1}$ maps $B^{k+1}_i$ to a subcomplex of
$B^{k+1}_{i+1}$. There is an induced map $\Phi_{k+2}$ from the mapping
cylinder of $\Phi_{k+1}\co B^{k+1}_i \to B^{k+1}_{i+1}$ to the mapping
cylinder of $\Phi_{k+1} \co B^{k+1}_{i+1} \to B^{k+1}_{i+2}$;  use 
$\Phi_{k+1}\times \id$ on $B^{k+1}_{i} \times I$ and 
$\Phi_{k+1}$ on $B^{k+1}_{i+1}$. Then  
$\Phi_{k+2}$ maps layer $i$ of $B^{k+2}_j$ to layer $i+1$ of
$B^{k+2}_{j+1}$ for any $i\leq j$ (in either hemisphere). These maps
defined on the layers of $B^{k+2}_j$ join together to define the map
$\Phi_{k+2} \co B^{k+2}_j \to B^{k+2}_{j+1}$. 

The balls $B^{k+2}_j$ embed into $\widetilde{X}^{k+2}_{r,P}$ exactly as
the balls $B^3_j$ embed into $\widetilde{X}^3_{r,P}$. That is, we
consider $\widetilde{X}^{k+2}_{r,P}$ as a union of copies of the mapping
cylinder of $\Phi_{k+1} \co \widetilde{X}^{k+1}_{r,P} \to
\widetilde{X}^{k+1}_{r,P}$ with the mapping parameter corresponding to
right multiplication by $u_k^{-1}$ or $v_k^{-1}$. Then the embedding
$B^{k+2}_j \to \widetilde{X}^{k+2}_{r,P}$ is assembled from the
embeddings $B^{k+1}_i \to \widetilde{X}^{k+1}_{r,P}$ (for $i \leq j$) as
shown in \fullref{fig:3ball}, with the upper hemisphere extending in
the $u_k$ direction and the lower hemisphere in the $v_k$
direction. Under this embedding, the map $\Phi_{k+2} \co B^{k+2}_j \to
B^{k+2}_{j+1}$ described above is simply the restriction of $\Phi_{k+2}
\co \widetilde{X}^{k+2}_{r,P} \to \widetilde{X}^{k+2}_{r,P}$ to
$B^{k+2}_j$. 

For any $k$, we define the \emph{shell} of a snowflake ball $B^k_j$ to be
the subspace $B^k_j - \Phi_{k}(B^k_{j-1})$, or simply $B^k_j$ in the case
$j=1$. 

\begin{lemma} \label{shellvolume} 
$\Vol^k(\shell(B^k_j)) \leq \Vol^{k-1}(\partial B^k_j)$. 
\end{lemma}

\begin{proof}
It suffices to show that every $k$--cell of the shell has a
$(k-1)$--dimensional face contained in $\partial B^k_j$. Recall that
$B^k_j$ is a union of layers, so consider the intersection of the shell
with layer $i$ (in either hemisphere). This layer is a mapping cylinder
$\mathcal{M}(\Phi_{k-1}\co B^{k-1}_i \to B^{k-1}_{i+1})$ and its
preimage in $B^k_{j-1}$ under $\Phi_{k}$ is layer $i-1$ of this
smaller ball (or is empty in the case $i=1$). Hence the intersection of
the shell with layer $i$ is 
\begin{equation*}
\begin{split} 
&\mathcal{M}(\Phi_{k-1}\co B^{k-1}_i \to B^{k-1}_{i+1}) \ - \ 
\Phi_{k}(\mathcal{M}(\Phi_{k-1} \co B^{k-1}_{i-1} \to
B^{k-1}_i)) \\
= \ & \mathcal{M}(\Phi_{k-1}\co B^{k-1}_i \to B^{k-1}_{i+1}) \ - \ 
\mathcal{M}(\Phi_{k} \co \Phi_{k-1}(B^{k-1}_{i-1}) \to
\Phi_{k-1}(B^{k-1}_i)))\\
= \ & \mathcal{M}(\Phi_{k-1}\co (B^{k-1}_i -
\Phi_{k-1}(B^{k-1}_{i-1})) \to (B^{k-1}_{i+1} -
\Phi_{k-1}(B^{k-1}_i))) \\
\end{split}
\end{equation*}
if $i>1$, and is $\mathcal{M}(\Phi_{k-1}\co B^{k-1}_i \to B^{k-1}_{i+1})$
in the case $i=1$. Either way, this part of $\shell(B^k_j)$ is the
mapping cylinder of the restriction of $\Phi_{k-1}$ to
$\shell(B^{k-1}_i)$. Hence each $k$--cell has a $(k-1)$--dimensional face
in $\shell(B^{k-1}_i)$, which is contained in $\partial B^k_j$. 
\end{proof} 

The next result is a higher-dimensional analogue of
\fullref{spherearea}. 

\begin{lemma} \label{boundaryvolumes}
Given $r$, $P$, and $k \geq 3$ there is a positive constant $F_k$ such
that $\Vol^{k-2}(\partial B^{k-1}_{j})$ 
$\leq \Vol^{k-1}(\partial B^k_j) \leq 
F_k \Vol^{k-2}(\partial B^{k-1}_{j})$ for every $j$. 
\end{lemma}
\begin{proof}
We prove, for $k \geq 3$, the following two statements: there exist
positive constants $E_k, F_k$ such that 
\begin{enumerate}
\item \label{s1} \ $(2C_1^{-1/\alpha}) (r^{1/\alpha})^j \ \leq
\ \Vol^{k-2}(\partial B_j^{k-1}) \ \leq \ E_k (r^{1/\alpha})^j$, and 
\item \label{s2} \  
$\Vol^{k-2} (\partial B_j^{k-1}) \ \leq \ \Vol^{k-1}
(\partial B_j^k) \ \leq \ F_k \Vol^{k-2}(\partial B_j^{k-1})$ 
\end{enumerate}
for all $j$ (with $C_1$ given by \fullref{sflength}). Statement
\eqref{s1} is a higher-dimensional analogue of 
\fullref{sflength} and \eqref{s2} is the main statement of the
lemma. The two statements are proved together by induction on $k$. 

If $k=3$ then \eqref{s1} follows from \fullref{sflength}, with
$E_3 = 2C_0^{-1/\alpha}$. Statement \eqref{s2} is given by
\fullref{spherearea} (with $F_3 = F_0$). 

For $k > 3$ we prove \eqref{s1} as follows. The induction hypothesis
implies that 
\[\Vol^{k-2}(\partial B_j^{k-1}) \ \leq \ F_{k-1} \Vol^{k-3}
(\partial B_j^{k-2})\]
by \eqref{s2} and $\Vol^{k-3}(\partial B_j^{k-2})
\leq E_{k-1} (r^{1/\alpha})^j$ by \eqref{s1}. Hence $\Vol^{k-2}(\partial
B_j^{k-1}) \leq E_k (r^{1/\alpha})^j$ with $E_k = F_{k-1} E_{k-1}$. We
also have (by induction) $\Vol^{k-2}(\partial B_j^{k-1}) \geq
\Vol^{k-3}(\partial B_j^{k-2}) \geq (2C_1^{-1/\alpha}) (r^{1/\alpha})^j$
by \eqref{s2} and \eqref{s1}. This establishes \eqref{s1}. 

To prove \eqref{s2} we count vertical and horizontal $(k-1)$--cells of
$\partial B^k_j$ as in the proof of \fullref{spherearea}. In each 
hemisphere of $B^k_j$, layer $i$ is a copy of the mapping cylinder of
$\Phi_{k-1} \co B^{k-1}_i \to B^{k-1}_{i+1}$. This layer meets $\partial
B^k_j$ in \emph{horizontal} cells which are the $(k-1)$--cells of
$\shell(B^{k-1}_i)$, and 
\emph{vertical} cells, each of which is the product of a $(k-2)$--cell in
$\partial B^{k-1}_i$ with $I$. This latter observation implies the first
inequality of \eqref{s2} (taking $i = j$) and also that the number of
vertical cells in layer $i$ is at most $\Vol^{k-2}(\partial
B^{k-1}_i)$. The number of horizontal cells is at most
$\Vol^{k-2}(\partial B^{k-1}_i)$ by \fullref{shellvolume}. Adding
the contributions from all layers in both hemispheres, we obtain 
\[\Vol^{k-1}(\partial B_j^k) \ \leq \ 4 \, \sum_{i=1}^j \Vol^{k-2}
(\partial B_i^{k-1}).\] 
Statement \eqref{s1} implies $4\sum_{i=1}^j \Vol^{k-2} (\partial
B_i^{k-1}) \leq 4 E_k\sum_{i=1}^j (r^{1/\alpha})^i$ and the latter
sum is a geometric series. Hence $\Vol^{k-1}(\partial B_j^k) \leq F'_k
(r^{1/\alpha})^j$ for some constant $F'_k$. Now \eqref{s1} implies that 
$\Vol^{k-1}(\partial B_j^k) \leq (F'_k/2)(C_1^{1/\alpha})
\Vol^{k-2}(\partial B_j^{k-1})$, establishing \eqref{s2} with $F_k =
(F'_k/2)C_1^{1/\alpha}$. 
\end{proof}

\section{Proof of Theorem C}\label{s:higher} 

We will establish upper and lower bounds for the $k$--dimensional Dehn
functions $\delta^{(k)}(x)$ of the groups $\Sigma^{k-1}G_{r,P}$ and these
will be equivalent. As usual $\lambda$ denotes the Perron--Frobenius 
eigenvalue of $P$ and $\alpha = \log_{\lambda}(r)$. In the case of
$\Sigma^{k-1}\Z^2$ we define $\alpha = 1$. 

\subsection*{The lower bound} 
As in the proof of  \fullref{mainthm1}, we show that the embedded
snowflake balls $B^{k+1}_i \subset \widetilde{X}^{k+1}_{r,P}$ have
the correct proportions and are numerous enough to determine 
$\delta^{(k)}(x)$ from below. 

First we show that for every $k\geq 1$ there is a constant $G_k$ such
that 
\begin{equation}\label{volbj} 
\Vol^{k+1}(B^{k+1}_i) \geq G_k \Vol^{k}(\partial B^{k+1}_i)^{2\alpha} 
\end{equation}
for all $i$. The case $k=1$ was proved in \eqref{sfdisklowerbound} with
$G_1 = (C_0)^2 4^{-\alpha}$. For $k >1$ we proceed by induction. Note
that $\Vol^{k+1}(B^{k+1}_i) \geq \Vol^{k}(B^{k}_i)$ since the latter is the
volume of the mapping cylinder of $\Phi_{k}\co B^{k}_i \to
\Phi_{k}(B^{k}_i)$ inside $B^{k+1}_i$. 
We also have $\Vol^{k}(B^{k}_i) \geq G_{k-1} \Vol^{k-1}(\partial
B^{k}_i)^{2\alpha}$ by the induction
hypothesis. \fullref{boundaryvolumes} implies that $G_{k-1}
\Vol^{k-1}(\partial B^{k}_i)^{2\alpha} \geq G_{k-1} F_{k+1}^{-2\alpha}
\Vol^{k}(\partial B^{k+1}_i)^{2\alpha}$. Equation \eqref{volbj} now follows
by taking $G_k = G_{k-1}F_{k+1}^{-2\alpha}$. 

Next we show that for each $k \geq 2$ the sequence $(\Vol^{k}(\partial
B^{k+1}_i))_i$ is exponentially bounded and tends to infinity. Consider first
the case $k=2$. Then we have
\[ \frac{\Vol^2(\partial B^3_{i+1})}{\Vol^2(\partial B^3_i)} \ \leq \ 
\frac{F_0\abs{\partial \Delta_{r^{i+1}}}}{\abs{\partial \Delta_{r^i}}}
\ \leq \ 
F_0 \Biggl( \frac{r^{i+1}C_1}{r^i C_0} \Biggr)^{1/\alpha} \ = \
\ F_0 \Biggl(\frac{rC_{1}}{C_{0}} \Biggr)^{1/\alpha} \]
where the first inequality holds by \fullref{spherearea}, and the second
by \fullref{sflength}. Thus, the sequence is exponentially bounded. For
$k>2$ we have 
\[ \frac{\Vol^{k}(\partial B^{k+1}_{i+1})}{\Vol^{k}(\partial B^{k+1}_i)}
\ \leq \ \frac{F_{k+1} \Vol^{k-1}(\partial
B^{k}_{i+1})}{\Vol^{k-1}(\partial B^{k}_i)}\]
by \fullref{boundaryvolumes} and so $(\Vol^{k}(\partial B^{k+1}_i))_i$ is
exponentially bounded, by induction on $k$. It tends to infinity because
\[\Vol^{k}(\partial B^{k+1}_i) \geq \Vol^2(\partial B^3_i) \geq
\abs{\partial \Delta_{r^i}} \geq 2C_1^{-1/\alpha} (r^{1/\alpha})^i\]
by \fullref{boundaryvolumes}, \fullref{spherearea}, and 
\fullref{sflength}. 
Now, using Remarks~\ref{sparse} and~\ref{k-embedded}, we conclude from
\eqref{volbj} that $\delta^{(k)}(x) \succeq x^{2\alpha}$. 

\subsection*{The upper bound} 
To establish the upper bound we must work with Dehn functions
$\delta^M_G(x)$ modeled on arbitrary manifolds $M$ with boundary, as
defined in \fullref{prelimsect}. Recall that the \emph{dimension} of
$\delta^M_G(x)$ is the dimension of $\partial M$, and $\delta^M_G(x)$ 
agrees with the usual $k$--dimensional Dehn function when $M$ is the
$(k+1)$--dimensional ball. 

A function $F \co \N \to \N$ is \emph{superadditive} if
$F(a+b) \geq F(a) + F(b)$ for all $a, b$. Recall that
the \emph{geometric dimension} of a group $G$ is the smallest dimension
of a $K(G,1)$ complex. 

\begin{theorem}\label{suspthm} 
Let $G$ be a group of type $\mathcal{F}_n$ and geometric dimension at
most $n$, and fix a finite aspherical $n$--complex $X$ with fundamental
group $G$. Suppose that the Dehn function $\delta^M_G(x)$ (defined with
respect to $X$) satisfies 
\[\delta^M_G(x) \leq F(x)\] 
for every $n$--manifold $M$, where $F\co \N \to \N$ is non-decreasing. Let
$H$ be a multiple ascending 
HNN extension of $G$. Then $H$ is of type $\mathcal{F}_{n+1}$, has
geometric dimension at most $n+1$, and 
\[\delta^M_H(x) \leq F(x)\]
for every $(n+1)$--manifold $M$. 
\end{theorem}
In the hypotheses we are including Dehn functions $\delta^M_G(x)$ where
$M$ has more than one connected component (otherwise we should add that
$F$ is superadditive). \textbf{Stipulation:} the $n$--dimensional Dehn
functions in the 
conclusion are defined with respect to a fixed complex $Y$ constructed in
the proof of the theorem. 

\begin{proof}
First we define the finite $(n+1)$--dimensional
complex $Y$ with fundamental group $H$ in the usual way. Suppose the
multiple ascending extension has $k$ stable letters. Form $k$ copies of
$X \times [-1,1]$, give each the product cell structure, and attach each copy
of $X \times \{-1\}$ to $X$ by the identity map. Then attach each copy of
$X \times \{1\}$ to $X$ by the appropriate monodromy map, and call
the resulting space $Y$. Let $Z \subset Y$ be the union of the spaces $X
\times \{0\}$. There are natural projections along the fibers $p_0
\co Z \to X$ and $p_1 \co Z \to X$ which factor through $Z \times
\{-1\}$ and $Z \times \{1\}$ respectively. Let $\widetilde{Y}$ be the
universal cover of $Y$ and let $\widetilde{X}$ and $\widetilde{Z}$ be the
preimages of $X$ and $Z$ in $\widetilde{Y}$. The projections $p_i$ lift
to projections $p_i \co \widetilde{Z} \to \widetilde{X}$ along
fibers. Note that each component of $\widetilde{X}$ and $\widetilde{Z}$
is a copy of the universal cover of $X$, and in fact $p_0 \co
\widetilde{Z} \to \widetilde{X}$ is a homeomorphism. 

Each open $k$--cell $\sigma^k$ in $\widetilde{Z} \times (-1,1) \subset
\widetilde{Y}$ has the form $\sigma^{k-1} \times (-1,1)$ where
$\sigma^{k-1}$ is a $(k-1)$--cell in $\widetilde{X}$, and the restriction
of $p_0$ to $\sigma^k \cap \widetilde{Z}$ is simply projection onto the
first factor. Since $\widetilde{Z}$ is not a subcomplex of
$\widetilde{Y}$, we measure volume in $\widetilde{Z}$ by passing to
$\widetilde{X}$ via $p_0$. The description of $p_0$ just given leads
to the following observation: if $f \co M^k \to \widetilde{Y}$ is an
admissible map transverse to $\widetilde{Z}$ and $\widetilde{X}$, and $N
= f^{-1}(\widetilde{Z})$ and $M_0 = f^{-1}(\widetilde{X})$, then 
$p_0 \circ f\vert_N$ and $f\vert_{M_0}$ are admissible and 
\begin{equation}\label{volume} 
\Vol^k(f) = \Vol^{k-1}(p_0 \circ f\vert_N) + \Vol^k(f \vert_{M_0}) 
\end{equation}
where the left hand side is volume in $\widetilde{Y}$ and the right hand
side is volume in $\widetilde{X}$. 

Now suppose that $M$ is a compact $(n+1)$--manifold with
boundary and let $g \co M \to \widetilde{Y}$ be a least-volume map with
boundary $f = g \vert_{\partial M}$. 
We can arrange by a homotopy that $N = g^{-1}(\widetilde{Z})$ is a
properly embedded codimension one submanifold with a product neighborhood
$N \times [-1,1] \subset M$ such that $g^{-1}(\widetilde{Z} \times
(-1,1)) = N \times (-1,1)$. (To do this, consider the composition
$\pi \co \widetilde{Y} \to Y \to S^1 \vee S^1 \to S^1$, where $S^1 \vee
S^1$ is the underlying graph for the multiple HNN-description of
$\pi_1(Y)$ and $S^1 \vee S^1 \to S^1$ is a fold. By a homotopy of $g$,
$\pi \circ g$ can be made smooth in a neighborhood of
$g^{-1}(\widetilde{Z})$. By a further homotopy, we can arrange that
$\pi(\widetilde{Z})$ is a regular value of $\pi \circ g$; now use
transversality.) 
The product structure on $N \times [-1,1]$ may be chosen so that
$g \vert_{N \times (-1,1)}$ is the map $g\vert_N \times \id$. Note that
$N$ may have several connected components. 

We claim that $\Vol^n(p_0 \circ g \vert_N)$ is smallest among all
$N$--fillings of $p_0 \circ f \vert_{\partial N} \co \partial N \to
\widetilde{X}$. Assuming this for the moment, the theorem is proved as
follows. We have $\Vol^{n+1}(g) = \Vol^n(p_0 \circ g\vert_N)$ by
\eqref{volume} because $\widetilde{X}$ has dimension $n$. Then
$\Vol^n(p_0 \circ g\vert_N) = \FVol^N(p_0 \circ f\vert_{\partial N})$
by \eqref{fvoleqn} and the claim, and the latter is at most
$\delta^N_G(\Vol^{n-1}(p_0 
\circ f\vert_{\partial N}))$ by the definition of $\delta^N_G$. Equation
\eqref{volume} implies that $\delta^N_G(\Vol^{n-1}(p_0 \circ
f\vert_{\partial N})) \leq \delta^N_G(\Vol^n(f))$. Then we have the
desired bound 
\[ \FVol^M(f) \ = \ \Vol^{n+1}(g) \ \leq \ \delta^N_G(\Vol^n(f)) \ \leq 
\ F(\Vol^n(f))\] 
by the main hypothesis and we conclude that $\delta^M_H(\Vol^n(f)) \leq 
F (\Vol^n(f))$. Since $\Vol^n(f)$ was arbitrary and $F$ is
non-decreasing, we have $\delta^M_H(x) \leq F(x)$ for all $x$. 

Now we return to the claim that $\Vol^n(p_0 \circ g\vert_N) =
\FVol^N(p_0 \circ f\vert_{\partial N})$. We show that if $p_0 \circ
g\vert_N$ is not a least-volume filling of $p_0 \circ f\vert_{\partial
N}$ then $g$ can be modified rel $\partial M$ to a map of smaller
volume, contradicting the choice of $g$. 

Let $M_0 = g^{-1}(\widetilde{X})$, and note that the frontier of $M_0$ in
$M$ is $N \times \{-1\} \cup N \times \{1\}$. These two subsets of
$\partial M_0$ will be denoted $M_0^-$ and $M_0^+$ respectively.  

Suppose $\Vol^n(h) < \Vol^n(p_0 \circ g \vert_N)$ for some map $h \co N
\to \widetilde{X}$ with $h \vert_{\partial N} = p_0 \circ f
\vert_{\partial N}$. Form a new copy of $M$ in which $N \times (-1,1)$ is
replaced by $N \times (-2,2)$. Define a new map $g'\co M \to
\widetilde{Y}$ by letting $g'$ be $g$ on $M_0$, $(p_0^{-1} \circ h)
\times \id$ on $N \times (-1,1)$, and by extending to the remaining
regions as follows. Note that $(p_0^{-1} \circ h) \times \id$ extends
continuously to $N \times [-1,1]$ as $h$ on $N \times \{-1\}$ and as
$p_1 \circ p_0^{-1} \circ h$ on $N \times \{1\}$. 
Since each component of $\widetilde{X}$ is contractible the maps $p_1
\circ p_0^{-1} \circ h$ and $g\vert_{M_0^+}$ are homotopic rel
$\partial N$. We let $g' \vert_{N \times [1,2]} \co N \times [1,2] \to
\widetilde{X}$ be such a homotopy. Similarly $g' \vert_{N \times [-2,
-1]}$ is defined to be a homotopy in $\widetilde{X}$ from $g
\vert_{M_0^-}$ to $h$, fixing $\partial N$ pointwise. This defines the
map $g' \co M \to \widetilde{Y}$. 

Now collapse each fiber of $\partial N \times [1,2]$ and $\partial N
\times [-2,-1]$ to a point, to obtain a new copy of $M$ with a map
$g''\co M \to \widetilde{Y}$ which agrees with $g$ on $\partial M$. 
Note that all of $M - (N \times (-1,1))$ maps by $g''$ into
$\widetilde{X}$ and $g'' \vert_{N \times (-1,1)} = (p_0^{-1} \circ
h) \times \id$. Hence by \eqref{volume} we have $\Vol^{n+1}(g'') =
\Vol^n(h) < \Vol^n(p_0 \circ g\vert_N) = \Vol^{n+1}(g)$, a
contradiction. 
\end{proof}

\begin{lemma}\label{2dim} 
If $G$ is finitely presented, $\delta_G(x) \leq F(x)$ with $F(x)$
superadditive, and $M$ is a compact $2$--manifold with
boundary, then $\delta^M_G(x) \leq F(x)$. 
\end{lemma}

In particular if $\delta_G(x)$ is superadditive then $\delta^M_G(x) \leq
\delta_G(x)$ for every compact $2$--manifold $M$. 

\begin{proof}
If $M$ is connected with one boundary component then let $q\co M \to
D^2$ be a quotient map which collapses the complement of a collar
neighborhood of $\partial M$ to a point. Then $\area(g \circ q) =
\area(g)$ for any map $g \co D^2 \to \widetilde{X}$, and we have
$\delta^M_G(x) \leq \delta_G(x) \leq F(x)$. 

If $N$ is closed then $\delta^{M \sqcup N}_G(x) = \delta^M_G(x)$ since
$N$ may be assigned zero area by mapping it to a point. So without
loss of generality assume that $M$ has no closed components. For each
component $M'$ of $M$ there is a quotient map to a connected, simply
connected space $Z'$ which is a union of disks (one for each
boundary component of $M'$) and arcs joining them. Taking a union of such
spaces and maps, we have a quotient map $M \to Z$. 
Every map $D^2 \sqcup \cdots \sqcup D^2 \to \widetilde{X}$ extends to a
map $Z \to \widetilde{X}$ which yields (by composition) a map $M\to
\widetilde{X}$ with the same area. Hence 
$\delta^M_G(x) \leq \delta^{D^2 \sqcup \cdots \sqcup D^2}_G(x)$. Now 
superadditivity of $F$ implies $\delta^{D^2 \sqcup \cdots \sqcup
D^2}_G(x) \leq F(x)$. 
\end{proof}

\begin{theorem}\label{suspcor} 
Let $G$ be a finitely presented group of geometric dimension $2$ with
$\delta_G(x)$ equivalent to a superadditive function. Let $H$ be obtained
from $G$ by performing $n$ iterated multiple ascending HNN 
extensions. Then $\delta^{(n+1)}_H(x) \preceq \delta_G(x)$. 
\end{theorem}

The upper bound of \fullref{mainthm2} follows immediately, by
\fullref{mainthm1}.  

\begin{proof}
Let $F_0(x)$ be superadditive where $F_0(x) \simeq \delta_G(x)$. Then
$\delta_G(x) \leq F(x) = C F_0(Cx) + Cx$ for some $C$ and $F(x)$ is
superadditive. The result now follows directly from \fullref{2dim} and
\fullref{suspthm}. 
\end{proof}

The case $n=1$ of \fullref{suspcor} was proved by Wang and Pride
\cite{WP}, using a more direct method.

\section{Products with $\Z$} 

In this section we determine higher Dehn functions of $G \times \Z$
for certain groups $G$. In these cases the geometry of $G\times \Z$ is
accurately represented by embedded balls which are products of optimal
balls in $G$ with intervals, with suitably chosen lengths. We conclude
the section by proving \fullref{mainthm3}. 

To establish an upper bound for Dehn functions of $G \times \Z$ we need
the following refinement of \fullref{suspthm}. The proof
is based on Theorem 6.1 of Alonso {\em et al} \cite{alonso+}. 

\begin{theorem}\label{upperGxZ} 
Let $G$ be a group of type $\mathcal{F}_n$ and geometric dimension at
most $n$, and fix a finite aspherical $n$--complex $X$ with fundamental group
$G$. Suppose that the Dehn function $\delta^M_G(x)$ satisfies
\[ \delta^M_G(x) \ \leq \ C x^s \]
for every $n$--manifold $M$, and fixed $C>0$ and $s > 1$. Then 
\[ \delta^M_{G \times \Z}(x) \ \leq \ C^{1/s} x^{2-1/s}\]
for every $(n+1)$--manifold $M$. 
\end{theorem}

\begin{proof} 
First note that we are in the situation of \fullref{suspthm}, which
is valid, but no longer provides the best possible upper bound. Define $Y$, 
$Z$, $p_0$, and $p_1$ as in the proof of \fullref{suspthm}. Note
that now the projections along fibers $p_0$, $p_1 \co \widetilde{Z}
\to \widetilde{X}$ are both homeomorphisms, and $\Vol^k(p_0 \circ f)
\ = \ \Vol^k(p_1 \circ f)$ for any $f \co N^k \to \widetilde{Z}$. 

Given a compact $(n+1)$--manifold $M$ with boundary, consider a
map $f \co \partial M \to \widetilde{Y}$. Arrange that $L =
f^{-1}(\widetilde{Z})$ is a codimension one submanifold with a product
neighborhood $L \times [-1,1] \subset \partial M$ such that
$f^{-1}(\widetilde{Z} \times (-1,1)) = L \times (-1,1)$. As before, the
product structure on $L \times [-1,1]$ can be chosen so that $f \vert_{L
\times (-1,1)}$ is the map $f\vert_L \times \id$. 

We will prove that $\delta^M_{G \times \Z}(x) \leq C^{1/s} x^{2-1/s}$ by
induction on the number of connected components of $L$. If $L =
\emptyset$ then $f(\partial M) \subset \widetilde{X}$. 
The components of $\partial M$ may map into different components of
$\widetilde{X}$. However, by joining these components with a minimal
collection of embedded arcs in the $1$--skeleton of $\widetilde{Y}$, one
obtains a contractible subcomplex $T\subset \widetilde{Y}$ of dimension
$n$ containing $f(\partial M)$. Then $f$ extends to a map $g \co M \to T
\subset \widetilde{Y}$ with $\Vol^{n+1}(g) = 0$. 

Now assume that $L \not= \emptyset$. 
Let $\widetilde{Z}_0$ be a connected component of $\widetilde{Z}$ such
that $L_0 = f^{-1}(\widetilde{Z}_0)$ is a non-empty union of components
of $L$, and $f(L)$ lies entirely in one component of $\widetilde{Y} - 
p_1(\widetilde{Z}_0)$. (Think of $L_0$ as an \emph{innermost} union of
components of $L$.) Let $N_1 \subset \partial M - (L_0 \times
(-1,1))$ be the union of components having boundary $L_0 \times
\{1\}$. That is, $N_1$ and its complement $N_{-1}$ in $\partial M - (L_0
\times (-1,1))$ map to opposite sides of $\widetilde{Z}_0 \times (-1,1)$
in $\widetilde{Y}$, and in fact $f(N_1) \subset p_1(\widetilde{Z}_0)
\subset \widetilde{X}$, by the choice of $\widetilde{Z}_0$. 

Our method now is to fill $L_0$ with a least-volume copy of $N_1$ and
then fill the two sides of $\partial M$ efficiently by $M$ (using
the induction hypothesis) and $N_1 \times I$. These fillings fit together
to yield a filling of $f$ by $M$ having the required volume. 

Let $v = \Vol^n(f)$ and $u = \Vol^{n-1}(p_0 \circ f\vert_{L_0})$ (which
is equal to $\Vol^n(f \vert_{L_0 \times (-1,1)})$ by \eqref{volume}). 
Let $h \co N_1 \to \widetilde{X}$ be a least-volume $N_1$--filling of
$p_0 \circ f\vert_{L_0}$. Thus, $h\vert_{\partial N_1} = p_0
\circ f\vert_{L_0}$ and $\Vol^n(h) \leq C u^s$. Define a new map $f' \co
\partial M \to  \widetilde{Y}$ by first collapsing the fibers of $L_0
\times [-1,1]$ to points, and then sending $N_{-1}$ by $f$ and $N_1$ by
$h$. Since $h$ is least-volume and $L_0 \times [-1,1]$ was collapsed we
have $\Vol^n(f') \leq v - u$. Also $(f')^{-1}(\widetilde{Z}) = L - L_0$, 
so by the induction hypothesis there is a map
$g_{-1} \co M \to \widetilde{Y}$ with $g_{-1} \vert_{\partial M} = f'$
such that 
\[ \Vol^{n+1}(g_{-1}) \leq C^{1/s}(v-u)^{2-1/s}.\] 
Next let $g_1 \co N_1 \times [-1,1] \to \widetilde{Y}$ be a homotopy 
which begins with $h$ on $N_1 \times \{-1\}$ and pushes across
$\widetilde{Z}_0 \times (-1,1)$ and then deforms within $p_1
(\widetilde{Z}_0)$ to $f \vert_{N_1}$, with the boundary fixed
pointwise. This latter homotopy exists since $p_1(\widetilde{Z}_0)$ is
contractible. Note that $\Vol^{n+1}(g_1) = \Vol^n(h)$ by \eqref{volume}
since $p_1(\widetilde{Z}_0)$ has dimension $n$. 

Now join $N_1 \subset \partial M$ to $(N_1 \times \{-1\}) \subset N_1
\times [-1,1]$ to get a new copy of $M$ and a map $g \co M \to
\widetilde{Y}$ extending $g_{-1}$ and $g_1$. Then $g \vert_{\partial M}
= f$ and 
\[ \Vol^{n+1}(g) \ \leq \ C^{1/s} (v-u)^{2-1/s} + v_h\] 
where $v_h = \Vol^n(h)$. Now $s > 1$ and $v \geq u$ imply 
\begin{equation}\label{eqn} 
\begin{split}
\Vol^{n+1}(g) \ &\leq \ C^{1/s}(v-u)v^{1-1/s} + v_h  \\
& = \ C^{1/s} v^{2-1/s}\left(1 - \frac{u}{v} +
\frac{v^{(1/s)-1}v_h}{C^{1/s} v}\right). 
\end{split}
\end{equation}
Recall that $v_h = \Vol^n(h) \leq \Vol^n(f \vert_{N_1}) \leq v$ because
$h$ is least-volume. Hence 
\begin{equation}\label{eqn2}
\begin{split}
1 - \frac{u}{v} + \frac{v^{(1/s)-1}v_h}{C^{1/s} v} \ 
&\leq \ 1 - \frac{u}{v} + \frac{v_h^{(1/s)-1}v_h}{C^{1/s} v} \\
&= \ 1 - \frac{u}{v} + \frac{v_h^{1/s}}{C^{1/s} v} \, . 
\end{split}
\end{equation}
The main hypothesis implies that $v_h \leq C u^s$, or
$v_h^{1/s} \leq C^{1/s}u$, again because $h$ is least-volume. Thus
\begin{equation}\label{eqn3} 
1 - \frac{u}{v} + \frac{v_h^{1/s}}{C^{1/s} v} \ \leq \ 1 - \frac{u}{v} +
\frac{u}{v} \ = \ 1. 
\end{equation}
By equations \eqref{eqn}, \eqref{eqn2}, and \eqref{eqn3} we have 
$\Vol^{n+1}(g) \leq C^{1/s}v^{2-1/s}$ where $v = \Vol^n(g\vert_{\partial
M})$, which completes the proof. 
\end{proof}

\begin{definition}\label{embeddedreps} 
Let $G$ be a group of type $\mathcal{F}_{k+1}$ and geometric dimension at
most $k+1$. The $k$--dimensional Dehn function $\delta^{(k)}_G(x)$
\emph{has embedded representatives} if there is a finite aspherical
$(k+1)$--complex $X$, a sequence of embedded $(k+1)$--dimensional balls
$B_i \subset \widetilde{X}$, and a function $F(x) \simeq
\delta^{(k)}_G(x)$, such that the sequence given by $(n_i) =
(\Vol^k(\partial B_i))$ tends to infinity and is exponentially bounded,
and $\Vol^{k+1}(B_i) \geq F(n_i)$ for each $i$. 
\end{definition}

The lower bounds established in this article for various Dehn functions
are all obtained by constructing embedded representatives and
applying Remarks~\ref{sparse} and \ref{k-embedded}. In particular the
$k$--dimensional Dehn functions of $\Sigma^{k-1}G_{r,P}$ and
$\Sigma^{k-1}\Z^2$ have embedded representatives. 

The next result generalizes \cite[Theorem 6.3]{alonso+} to higher
dimensions. 

\begin{proposition}\label{lowerGxZ} 
Let $G$ be a group of type $\mathcal{F}_{k+1}$ and geometric dimension at
most 
$k+1$. Suppose the $k$--dimensional Dehn function $\delta^{(k)}(x)$ of $G$
is equivalent to $x^s$ and has embedded representatives. Then $G \times
\Z$ has $(k+1)$--dimensional Dehn function $\delta^{(k+1)}(x)\succeq x^{2
- 1/s}$, with embedded representatives. 
\end{proposition}

\begin{proof}
We establish the lower bound $\delta^{(k+1)}(x) \succeq 
x^{2-1/s}$ for $G \times \Z$ as follows. Since
$\delta^{(k)}_G(x)$ has embedded representatives, let $X$, $F(x)$, $B_i$,
and  $(n_i)$ be as in \fullref{embeddedreps}; without loss of
generality suppose that $F(x) = Cx^s$ for some $C > 0$. Define $m_i = 3
\Vol^{k+1}(B_i)$. The space $Y = X \times S^1$ has fundamental group $G
\times \Z$ and universal cover $\widetilde{Y} = 
\widetilde{X} \times \R$. Consider the $(k+2)$--dimensional balls 
\[ C_i \ = \ B_i \times [0,m_i/3n_i] \ \subset \ \widetilde{Y}. \] 
The boundary of $C_i$ is $\partial B_i \times [0,m_i/3n_i] \cup B_i
\times \partial [0, m_i/3n_i]$ which implies that 
\[\Vol^{k+1}(\partial C_i) = m_i.\]
We also have $\Vol^{k+2}(C_i) = \Vol^{k+1}(B_i) m_i/3n_i =
(m_i)^2/9n_i$ for each $i$. Since $m_i = 3 \Vol^{k+1}(B_i) \geq
3C(n_i)^s$ we have $(3C)^{-1/s} (m_i)^{1/s} \geq n_i$. Then
\begin{equation*}
\Vol^{k+2}(C_i) \ = \ \frac{(m_i)^2}{9n_i} \ \geq
\ \left(\frac{C^{1/s}}{3^{2-1/s}}\right) (m_i)^{2-1/s}.
\end{equation*}
Note that $\widetilde{Y}$ is aspherical and has dimension $k+2$, and so
$C_i$ is a least-volume ball (cf \fullref{k-embedded}). Therefore
$\delta^{(k+1)}(m_i) \geq (C^{1/s} / 3^{2-1/s}) (m_i)^{2-1/s}$ for each
$i$. Now it remains to check that the sequence $(m_i)$ has the required
properties. It tends to infinity since $m_i \geq 3C(n_i)^s$. Also each
ball $B_i \subset \widetilde{X}$ is least-volume, so there is a 
constant $D$ such that $m_i \leq D(n_i)^s$ for all $i$.\footnote{Here we
are using the upper bound for $\delta^{(k)}_G(x)$.} Then $m_{i+1}/m_i
\leq (D/C)(n_{i+1}/n_i)^s$, which is bounded. Now \fullref{sparse}
implies that $\delta^{(k+1)}(x) \succeq x^{2-1/s}$. 
\end{proof}

We are now in a position to prove \fullref{mainthm3}. 

\begin{proof}[Proof of \fullref{mainthm3}] 
Fix $r$, $P$, and $q$, let $s(\ell) = \frac{2(\ell + 1)\alpha -
\ell}{2\ell \alpha - (\ell -1)}$, and let $G_{\ell}$ be the group
$\Sigma^{q-1}G_{r,P} \times \Z^{\ell}$. (Or let $s(\ell) = \frac{\ell +
2}{\ell + 1}$ and $G_{\ell} = \Sigma^{q-1}\Z^2 \times \Z^{\ell}$.) 
We verify by induction on $\ell$
the following statements for $G_{\ell}$: 
\begin{enumerate}
\item \label{st1} $\delta^M(x) \leq Cx^{s(\ell)}$ for all
$(q+\ell+1)$--manifolds $M$ and some constant $C > 0$, 
\item \label{st2} $\delta^{(q+\ell)}(x) \succeq x^{s(\ell)}$, and 
\item \label{st3} $\delta^{(q+\ell)}(x)$ has embedded representatives. 
\end{enumerate}
The first two statements together imply $\delta^{(q+\ell)}(x) \simeq
x^{s(\ell)}$. 

If $\ell = 0$ then \eqref{st1} follows from \fullref{suspthm} and
\fullref{2dim}. Statement~\eqref{st2} holds by
\fullref{mainthm2}, and we have already observed that \eqref{st3}
holds for these groups. 

For $\ell > 0$ note first that $s(\ell) = 2 - 1/s(\ell -1)$. Then
statement~\eqref{st1} holds by \fullref{upperGxZ} and 
property \eqref{st1} of $G_{\ell -1}$. \fullref{lowerGxZ} implies
\eqref{st2} and~\eqref{st3} by properties~\eqref{st1}--\eqref{st3} of
$G_{\ell -1}$. 
\end{proof}

\bibliographystyle{siam}
\bibliography{martinsbib} 

\def\cprime{$'$}
\begin{thebibliography}{10}

\bibitem{alonso+}
{\sc J.~M. Alonso, W.~A. Bogley, R.~M. Burton, S.~J. Pride, and X.~Wang}, {\em
  Second order {D}ehn functions of groups}, Quart. J. Math. Oxford Ser. (2), 49
  (1998), pp.~1--30.

\bibitem{AWP}
{\sc J.~M. Alonso, X.~Wang, and S.~J. Pride}, {\em Higher-dimensional
  isoperimetric (or {D}ehn) functions of groups}, J. Group Theory, 2 (1999),
  pp.~81--112.

\bibitem{BrBr}
{\sc N.~Brady and M.~R. Bridson}, {\em There is only one gap in the
  isoperimetric spectrum}, Geom. Funct. Anal., 10 (2000), pp.~1053--1070.

\bibitem{bf}
{\sc N.~Brady and M.~Forester}, {\em Density of isoperimetric spectra}.
\newblock Preprint, 2008.

\bibitem{mrb:bfs}
{\sc M.~R. Bridson}, {\em The geometry of the word problem}, in Invitations to
  geometry and topology, vol.~7 of Oxf. Grad. Texts Math., Oxford Univ. Press,
  Oxford, 2002, pp.~29--91.

\bibitem{mb:plms}
\leavevmode\vrule height 2pt depth -1.6pt width 23pt, {\em Polynomial {D}ehn
  functions and the length of asynchronously automatic structures}, Proc.
  London Math. Soc. (3), 85 (2002), pp.~441--466.

\bibitem{BuRoSa}
{\sc S.~Buoncristiano, C.~P. Rourke, and B.~J. Sanderson}, {\em A geometric
  approach to homology theory}, Cambridge University Press, Cambridge, 1976.
\newblock London Mathematical Society Lecture Note Series, No. 18.

\bibitem{burillo}
{\sc J.~Burillo}, {\em Lower bounds of isoperimetric functions for nilpotent
  groups}, in Geometric and computational perspectives on infinite groups
  (Minneapolis, MN and New Brunswick, NJ, 1994), vol.~25 of DIMACS Ser.
  Discrete Math. Theoret. Comput. Sci., Amer. Math. Soc., Providence, RI, 1996,
  pp.~1--8.

\bibitem{csc}
{\sc T.~Coulhon and L.~Saloff-Coste}, {\em Isop\'erim\'etrie pour les groupes
  et les vari\'et\'es}, Rev. Mat. Iberoamericana, 9 (1993), pp.~293--314.

\bibitem{ECHLPT}
{\sc D.~B.~A. Epstein, J.~W. Cannon, D.~F. Holt, S.~V.~F. Levy, M.~S. Paterson,
  and W.~P. Thurston}, {\em Word processing in groups}, Jones and Bartlett
  Publishers, Boston, MA, 1992.

\bibitem{gromov1}
{\sc M.~Gromov}, {\em Hyperbolic groups}, in Essays in group theory, vol.~8 of
  Math. Sci. Res. Inst. Publ., Springer, New York, 1987, pp.~75--263.

\bibitem{gromov2}
\leavevmode\vrule height 2pt depth -1.6pt width 23pt, {\em Asymptotic
  invariants of infinite groups}, in Geometric group theory, Vol.\ 2 (Sussex,
  1991), vol.~182 of London Math. Soc. Lecture Note Ser., Cambridge Univ.
  Press, Cambridge, 1993, pp.~1--295.

\bibitem{HV}
{\sc A.~Hatcher and K.~Vogtmann}, {\em Isoperimetric inequalities for
  automorphism groups of free groups}, Pacific J. Math., 173 (1996),
  pp.~425--441.

\bibitem{KH}
{\sc A.~Katok and B.~Hasselblatt}, {\em Introduction to the modern theory of
  dynamical systems}, vol.~54 of Encyclopedia of Mathematics and its
  Applications, Cambridge University Press, Cambridge, 1995.
\newblock With a supplementary chapter by Katok and Leonardo Mendoza.

\bibitem{ol}
{\sc A.~Y. Ol{\cprime}shanskii}, {\em Groups with quadratic-non-quadratic
  {D}ehn functions}, Internat. J. Algebra Comput., 17 (2007), pp.~401--419.

\bibitem{olsap}
{\sc A.~Y. Ol{\cprime}shanskii and M.~V. Sapir}, {\em Groups with small {D}ehn
  functions and bipartite chord diagrams}, Geom. Funct. Anal., 16 (2006),
  pp.~1324--1376.

\bibitem{papa}
{\sc P.~Papasoglu}, {\em Isodiametric and isoperimetric inequalities for
  complexes and groups}, J. London Math. Soc. (2), 62 (2000), pp.~97--106.

\bibitem{BRS}
{\sc M.~V. Sapir, J.-C. Birget, and E.~Rips}, {\em Isoperimetric and
  isodiametric functions of groups}, Ann. of Math. (2), 156 (2002),
  pp.~345--466.

\bibitem{seneta}
{\sc E.~Seneta}, {\em Nonnegative matrices and {M}arkov chains}, Springer
  Series in Statistics, Springer-Verlag, New York, second~ed., 1981.

\bibitem{vscc}
{\sc N.~T. Varopoulos, L.~Saloff-Coste, and T.~Coulhon}, {\em Analysis and
  geometry on groups}, vol.~100 of Cambridge Tracts in Mathematics, Cambridge
  University Press, Cambridge, 1992.

\bibitem{wang}
{\sc X.~Wang}, {\em Second order {D}ehn functions of finitely presented groups
  and monoids}, 1996.
\newblock PhD Thesis, {U}niversity of {G}lasgow.

\bibitem{WP}
{\sc X.~Wang and S.~J. Pride}, {\em Second order {D}ehn functions and
  {HNN}-extensions}, J. Austral. Math. Soc. Ser. A, 67 (1999), pp.~272--288.

\bibitem{young}
{\sc R.~Young}, {\em A note on higher-order filling functions}.
\newblock Preprint, arxiv:0805.0584.

\end{thebibliography}

\end{document}